\documentclass[12pt]{article}
\usepackage{amsmath}
\usepackage{amsthm}
\usepackage{amsfonts}
\usepackage{amssymb}
\usepackage{graphicx}
\usepackage{xypic}
\usepackage{a4wide}
\usepackage{enumerate}
\usepackage{mathrsfs}
\newtheorem{theorem}{Theorem}[section]
\newtheorem{proposition}[theorem]{Proposition}
\newtheorem{corollary}[theorem]{Corollary}
\newtheorem{lemma}[theorem]{Lemma}
\theoremstyle{definition}

\newtheorem{definition}[theorem]{Definition}
\newtheorem{example}[theorem]{Example}

\newtheorem{remark}[theorem]{Remark}

\newcommand{\bcomod}[2]{{}^{#1}\mathcal{M}^{#2}}
\newcommand{\cat}[1]{\mathbf{#1}}
\newcommand{\cohom}[3]{\mathrm{h}_{#1}(#2,#3)}
\newcommand{\coring}[1]{\mathfrak{#1}}
\newcommand{\cotensor}[1]{\square_{#1}}
\renewcommand{\hom}[3]{\mathrm{Hom}_{#1}(#2,#3)}
\newcommand{\lcomod}[1]{{}^{#1}\mathcal{M}}
\newcommand{\rcomod}[1]{\mathcal{M}^{#1}}
\newcommand{\rmod}[1]{\mathcal{M}_{#1}}
\newcommand{\lmod}[1]{_{#1}\mathcal{M}}
\newcommand{\ldual}[1]{{ }^{\ast}#1}
\newcommand{\rdual}[1]{#1^{\ast}}
\newcommand{\lDual}{{ }^{\ast}(-)}
\newcommand{\rDual}{ (-)^{\ast}}
\newcommand{\tensor}[1]{\otimes_{#1}}
\newcommand{\e}[2]{\mathrm{e}_{#1}(#2)}
\newcommand{\h}[2]{\mathrm{H}(#1,#2)}

\begin{document}

\title{Comatrix Coring generalized and Equivalences of
Categories of Comodules}
\author{Mohssin Zarouali-Darkaoui\\
\normalsize Departamento de \'{A}lgebra \\
\normalsize Facultad de Ciencias \\
\normalsize Universidad de Granada\\
\normalsize E18071 Granada, Spain \\
\normalsize E-mail: \textsf{zaroual@correo.ugr.es}}

\date{\empty}

\maketitle

\section*{Introduction}

Corings were introduced by Sweedler in \cite{Sweedler:1975}. A
coring over an associative algebra with unit over a commutative ring
$k$, $A$, is an $A$-bimodule $\coring{C}$ with two $A$-bimodule maps
$\Delta_{\coring{C}}:\coring{C}\rightarrow\coring{C}\tensor{A}{\coring{C}}$
(comultiplication or coproduct) and
$\epsilon_{\coring{C}}:\coring{C}\rightarrow A$ (counit) such that
the same diagrams as for coalgebras are commutative. Recently,
corings were intensively studied. The main motivation of this
studies is the observation of Takeuchi,
\cite[32.6]{Brzezinski/Wisbauer:2003}, that relates entwining
structures (resp. entwined modules) with corings (resp. comodules
over corings). For a detailed study of corings, we refer to
\cite{Brzezinski/Wisbauer:2003}.

Comatrix corings were introduced for the first time in
\cite{ElKaoutit/Gomez:2003}. In \cite{Brzezinski/Gomez:2003}, the
authors have given another definition and some properties of them.

In this paper we extend comatrix coring to the case of quasi-finite
comodules, and also we extend some of its interesting properties
given in \cite{Brzezinski/Gomez:2003}. Suppose that $_A\coring{C}$
and $_B\coring{D}$ are flat. Let
$X\in\bcomod{\coring{C}}{\coring{D}}$ and
$\Lambda\in\bcomod{\coring{D}}{\coring{C}}$ and suppose that
$-\square _{\coring{C}}X$ is a left adjoint to
$-\square_{\coring{D}}\Lambda$. If $\coring{C}_A$ and $\coring{D}_B$
are flat and $_{\coring{C}}X$, $_{\coring{D}}\Lambda$ are coflat, or
$A$ and $B$ are von Neumann, or if $\coring{C}$ and $\coring{D}$ are
coseparable, then the $B$-bimodule $\Lambda\cotensor{\coring{C}}X$
is endowed by a structure of $B$-coring (Theorem
\ref{Comatrixcoring}). The coproduct is
\[
\xymatrix{\Delta:\Lambda\square_{\coring{C}}X \ar[r]^-\simeq &
\Lambda\square_{\coring{C}}(
\coring{C}\square_{\coring{C}}X)\ar[rr]^-{\Lambda\square
_{\coring{C}}(\mathcal{\psi}\square_{\coring{C}}X)} & &
\Lambda\square_{\coring{C}}((X\square_{\coring{D}}\Lambda)
\square_{\coring{C}}X)\ar[r]^-\simeq &
(\Lambda\square_{\coring{C}}X)\square_{\coring{D}}(
\Lambda\square_{\coring{C}}X) \ar@{^{(}->}[d] \\ & & & &
(\Lambda\square_{\coring{C}}X)\otimes_B(\Lambda
\square_{\coring{C}}X),}
 \]
 and the counit is
$\xymatrix{\epsilon=\epsilon_{\coring{D}}\circ\omega
:\Lambda\square_{\coring{C}}X\ar[r] & B}$.

The most important property of this coring is the following. Let
$\Lambda\in \bcomod{\coring{D}}{\coring{C}}$ be a bicomodule,
quasi-finite as a right $\coring{C}$-comodule, such that
$_A\coring{C}$ and $_B\coring{D}$ are flat. Set
$X=\cohom{\coring{C}}{\Lambda}{\coring{C}}\in
\bcomod{\coring{D}}{\coring{C}}$.
 If
 \begin{enumerate}[(a)]\item $\coring{C}_A$ and $\coring{D}_B$ are flat,
 the cohom functor $\cohom{\coring{C}}{\Lambda}{-}$ is exact and
$_{\coring{D}}{\Lambda}$ is coflat, or
\item $A$ and $B$ are von Neumann regular rings and the cohom functor
$\cohom{\coring{C}}{\Lambda}{-}$ is exact, or
\item $\coring{C}$ and
$\coring{D}$ are coseparable corings,
\end{enumerate}
then the canonical isomorphism
$\xymatrix{\delta_{\Lambda}:\e{\coring{C}}{\Lambda}\ar[r] &
\Lambda\cotensor{\coring{C}}X}$ is an isomorphism of $B$-corings
(see Proposition \ref{comatrixproperty1}). We think that, under
these conditions, this coring isomorphism give a concrete
description of the coendomorphism coring even in the case of
coalgebras over fields. Indeed, let $C$ be a coalgebra over a field
$k$, and let $\Lambda\in\rcomod{C}$ be a quasi-finite and injective
comodule. Then $\e{C}{\Lambda}\simeq
\Lambda\cotensor{C}\cohom{C}{\Lambda}{C}$ as coalgebras (the case
(b)). We also think that the last isomorphism give a more concrete
characterization of equivalence between categories of comodules over
coalgebras over a field, than that of Takeuchi \cite{Takeuchi:1977}.
of course, this characterization of the coendomorphism coring is
also useful in the characterization of equivalence between
categories of comodules over corings. Indeed, in Section 3, we study
equivalences between categories of comodules over rather general
corings. We generalize and improve (using results we give in
\cite{Zarouali:2004}) the main results concerning equivalences
between categories of comodules given in \cite{Takeuchi:1977},
\cite{AlTakhman:2002} and \cite{Brzezinski/Wisbauer:2003}. We also
give new characterizations of equivalences between categories of
comodules over coseparable corings or corings \cite{Zarouali:2004}
with a duality. We apply our results to the particular case of the
adjoint pair of functors associated to a morphism of corings over
different base rings \cite{Gomez:2002}. Finally, when applied to
corings associated to entwining structures, and that associated to a
$G$-graded algebra and a right $G$-set, we obtain new results
concerning entwined modules and graded modules. We think that our
result, Theorem \ref{graded2}, is more simple than that Del R\'{\i}o
\cite[Theorem 2.3]{DelRio:1992}.

The paper is organized as follows. In Section \ref{basic}, we give
some useful definitions and notations. In Section
\ref{comatrixgeneralized}, we generalize comatrix corings introduced
in \cite{ElKaoutit/Gomez:2003} and \cite{Brzezinski/Gomez:2003}, and
we generalize also some properties given in
\cite{Brzezinski/Gomez:2003}. Section \ref{equivalences} is devoted
to the study of equivalences of comodule categories over corings.
Our results given in \cite{Zarouali:2004} and Section
\ref{comatrixgeneralized}, will allow us to generalize and improve
the main results in both \cite{AlTakhman:2002} and
\cite{Brzezinski/Wisbauer:2003} (see Propositions
\ref{equivalence1}, \ref{equivalence2}, Theorems \ref{equivalence3},
\ref{equivalence4}), and also to give new results concerning
equivalences of comodule categories over coseparable corings (see
Propositions \ref{equivalence1}, \ref{equivalence2}, Theorems
\ref{equivalence3}, \ref{equivalence5}) and over corings over QF
rings (Theorem \ref{equivalence6}). Obviously our last theorem
generalizes \cite[Theorem 3.5]{Takeuchi:1977} and \cite[Corollary
7.6]{AlTakhman:2002}. In Section \ref{induction}, we deals with the
application of some of our results giving in Section
\ref{equivalences} to the induction functor. In Section
\ref{entwined}, we apply some of our results giving in the previous
sections to the corings associated to entwining structures, and in
particular those associated to a $G$-graded algebra and a right
$G$-set, where $G$ is a group.

\section{Preliminaries and basic
notations}\label{basic}

Throughout this paper and unless otherwise stated, $k$ denote a
commutative ring (with unit), $A,$ $A',$ $A'',$ and $B$ denote
associative and unitary algebras over $k$, and $\coring{C},$
$\coring{C}',$ $\coring{C}'',$ and $\coring{D}$ denote corings over
$A,$ $A',$ $A'',$ and $B$, respectively.

A category $\cat{C}$ is said to be $k$-\emph{category} (called
$k$-linear category in \cite[I.0.1.2]{Saavedra:1972}) if for every
$M$ and $N$ in $\cat{C}$, $\hom{\cat{C}}{M}{N}$ is a $k$-module, and
the composition is $k$-bilinear. An abelian category which is a
$k$-category is said to be $k$-\emph{abelian category}. A functor
between $k$-categories is said to be $k$-\emph{functor} or
$k$-\emph{linear functor} if it is $k$-linear on the $k$-modules of
morphisms. A functor between $k$-categories is said to be a
$k$-\emph{equivalence} if it is $k$-linear and an equivalence.

We recall from \cite{Sweedler:1975} that an $A$-\emph{coring}
consists of an $A$-bimodule $\coring{C}$ with two $A$-bimodule maps
\[
\Delta : \coring{C} \rightarrow \coring{C} \tensor{A} \coring{C},
\quad \epsilon : \coring{C} \rightarrow A
\]
such that $(\coring{C} \tensor{A} \Delta) \circ \Delta = (\Delta
\tensor{A} \coring{C}) \circ \Delta$ and $(\epsilon \tensor{A}
\coring{C}) \circ \Delta = ( \coring{C} \tensor{A} \epsilon) \circ
\Delta = 1_\coring{C}$. A \emph{right} $\coring{C}$-\emph{comodule}
is a pair $(M,\rho_M)$ consisting of a right $A$-module $M$ and an
$A$-linear map $\rho_M: M \rightarrow M \tensor{A} \coring{C}$
(coaction) satisfying $(M \tensor{A} \Delta) \circ \rho_M = (\rho_M
\tensor{A} \coring{C}) \circ \rho_M$, and $(M \tensor{A} \epsilon)
\circ \rho_M = 1_M$. A \emph{morphism} of right
$\coring{C}$-comodules $(M,\rho_M)$ and $(N,\rho_N)$ is a right
$A$-linear map $f: M \rightarrow N$ such that $(f \tensor{A}
\coring{C}) \circ \rho_M = \rho_N \circ f$. The $k$-module of all
such morphisms will be denoted by $\hom{\coring{C}}{M}{N}$. Right
$\coring{C}$-comodules together with their morphisms form the
$k$-category $\rcomod{\coring{C}}$. Coproducts and cokernels (and
then inductive limits) in $\rcomod{\coring{C}}$ exist and they
coincide respectively with coproducts and cokernels in the category
of right $A$-modules $\rmod{A}$. If ${}_A\coring{C}$ is flat, then
it is easy to show that the subcomodules of $\coring{C}^n$, $n\in
\mathbb{N}$, form a family of generators. Hence, if $_A\coring{C}$
is flat, then the category $\rcomod{\coring{C}}$ is a Grothendieck
category. The converse is not true in general (see \cite[Example
1.1]{ElKaoutit/Gomez/Lobillo:2001unp}). When $\coring{C}=A$, with
the trivial $A$-coring structure, $\rcomod{A}$ is the category of
right $A$-modules $\rmod{A}$.

Now assume that the $A'-A$-bimodule $M$ is also a left comodule over
an $A'$-coring $\coring{C}'$ with structure map $\lambda_M : M
\rightarrow \coring{C}' \tensor{A'} M$. Assume moreover that
$\rho_M$ is $A'$-linear, and $\lambda_M$ is $A$-linear. It is clear
that $\rho_M : M \rightarrow M \tensor{A} \coring{C}$ is a morphism
of left $\coring{C}'$-comodules if and only if $\lambda_M : M
\rightarrow \coring{C}' \tensor{A'} M$ is a morphism of right
$\coring{C}$-comodules. In this case, we say that $M$ is a
$\coring{C}'-\coring{C}$-bicomodule. A morphism of bicomodules is a
morphism of right and left comodules. Then we obtain a $k$-category
$^\coring{C'}\mathcal{M}^\coring{C}$. If in particular
$\coring{C}'=A', \coring{C}=A$, then
$^\coring{C'}\mathcal{M}^\coring{C}$ is the category of
$A'-A$-bimodules $_{A'}\mathcal{M}_A$.

A coring $\coring{C}$ is said to be \emph{coseparable}
\cite{Guzman:1989} if the comultiplication map $\Delta_{\coring{C}}$
is a section in the category $^\coring{C}\mathcal{M}^\coring{C}$.
Obviously the trivial $A$-coring $\coring{C} = A$ is coseparable. We
refer to \cite{Brzezinski/Wisbauer:2003} for the definition and
basic properties of the notions: coendomorphism coring, cosplit
coring, and Frobenius coring.

Let $Z$ be a left $A$-module and $f: X \rightarrow Y$ a morphism in
$\rmod{A}$. Following \cite[40.13]{Brzezinski/Wisbauer:2003} we say
that $f$ is \emph{$Z$-pure} when the functor $-\tensor{A} Z$
preserves the kernel of $f$. If $f$ is $Z$-pure for every $Z \in
{}\lmod{A}$ then we say simply that $f$ is \emph{pure} in
$\rmod{A}$.

A bicomodule $N \in \bcomod{\coring{C}}{\coring{D}}$ is said to be
\emph{quasi-finite} as a right $\coring{D}$-comodule if the functor
$- \tensor{A} N : \rmod{A} \rightarrow \rcomod{\coring{D}}$ has a
left adjoint $\cohom{\coring{D}}{N}{-} : \rcomod{\coring{D}}
\rightarrow \rmod{A}$, and we call it the \emph{cohom functor}. If
$\omega_{Y,N}:=\rho_Y\otimes_AN-Y\otimes_A\lambda_N$ is $\coring{D}
\tensor{B} \coring{D}$-pure for every right $\coring{C}$-comodule
$Y$ (e.g., ${}_B\coring{D}$ is flat or $\coring{C}$ is coseparable)
then $N_{\coring{D}}$ is quasi-finite if and only if $-
\cotensor{\coring{C}} N : \rcomod{\coring{C}} \rightarrow
\rcomod{\coring{D}}$ has a left adjoint, which is also denoted by
$\cohom{\coring{D}}{N}{-}$ \cite[Proposition 4.2]{Gomez:2002}. In
particular ($\coring{D}=B$), a $\coring{C}-B$-bimodule $N$ is
quasi-finite as a right $B$-(co)module if and only if $_AN$ is
finitely generated and projective. The cohom functor is
$-\tensor{B}\ldual{N}: \rmod{B} \rightarrow \rmod{A}$, where
$\ldual{N}=\hom{A}{N}{A}$.

Recall from \cite[Exercise 19.19]{Anderson/Fuller:1992}, that a
module ${}_AW$ is said to be \emph{completely faithful} if
$$Ann_M(W):=\{m\in M\mid m\tensor{A}w=0 \; \textrm{in}
\; M\tensor{A}W \;\textrm{for all}\; w\in W\}=0,$$ for every right
$A$-module $M$. It follows from \cite[Proposition
II.7.2]{Mitchell:1965}, that a module ${}_AW$ is completely faithful
if and only if the functor $-\tensor{A}W$ is faithful.

Finally, the notation $\otimes$ will stand for the tensor product
over $k$.

\section{Comatrix coring generalized}\label{comatrixgeneralized}

In this Section we generalize the concept of comatrix coring defined
in \cite{Brzezinski/Gomez:2003} and \cite{ElKaoutit/Gomez:2003} to
the case of a quasi-finite comodule, and generalize some of its
properties.

\medskip
At first we will recall the definition of the \emph{cotensor
product} of comodules. Let
 $M\in{}\bcomod{\coring{C'}}{\coring{C}}$ and
$N\in{}^{\coring{C}}\mathcal{M}^{\coring{C}''}$. The map
\[
\omega_{M,N}=\rho_M\otimes_AN-M\otimes_A\lambda_N:M\otimes_AN
\rightarrow M\otimes_A\coring{C}\otimes_AN
\]
is a $\coring{C}'-\coring{C}^{\prime\prime}$-bicomodule map. Its
kernel in $_{A'}\mathcal{M}_{A''}$ is the \emph{cotensor product} of
$M$ and $N$, and it is denoted by $M\square_{\coring{C}}N$ . If
$\omega_{M,N}$ is $\coring{C}'_{A'}$-pure and
$_{A''}\coring{C}''$-pure, and the following
\begin{equation}
\operatorname{ker} (\omega_{M,N})\otimes_{A''}\coring{C}''
\otimes_{A''}\coring{C}'', \quad \coring{C}' \otimes_{A'}\coring{C}'
\otimes_{A'}\operatorname{ker} (\omega_{M,N}) \quad \textrm{and}
\quad \coring{C}'\otimes_{A'}\operatorname{ker}
(\omega_{M,N})\otimes_{A''}\coring{C}''
\end{equation}
are injective maps, then $M\square_{\coring{C}}N$ is the kernel of
$\omega_{M,N}$ in ${}^{\coring{C}'}\mathcal{M}^{\coring{C}''}$. This
is the case if $\omega_{M,N}$ is
$(\coring{C}'\otimes_{A'}\coring{C}')_{A'}$-pure,
$_{A''}(\coring{C}'' \otimes_{A''}\coring{C}'')$-pure, and
$\coring{C}'\otimes_{A'}\omega_{M,N}$ is $_{A''}\coring{C}''$-pure
(e.g. if $\coring{C}'_{A'}$ and $_{A''}\coring{C}''$ are flat, or if
$\coring{C}$ is a coseparable $A$-coring).

If for every $M\in{}\bcomod{\coring{C'}}{\coring{C}}$ and
$N\in{}^{\coring{C}}\mathcal{M}^{\coring{C}''}$, $\omega_{M,N}$ is
$\coring{C}'_{A'}$-pure and $_{A''}\coring{C}''$-pure, then we have
a $k$-linear bifunctor
\begin{equation}\label{cotensorbifunctor}
\xymatrix{-\cotensor{\coring{C}}-:{}^{\coring{C}'}\mathcal{M}^\coring{C}
\times{}^{\coring{C}}\mathcal{M}^{\coring{C}''}\ar[r]&
{}^{\coring{C}'}\mathcal{M}^{\coring{C}''}}.
\end{equation} If in particular $\coring{C}'_{A'}$ and
$_{A''}\coring{C}'' $ are flat, or if $\coring{C}$ is a coseparable
$A$-coring, then the bifunctor \eqref{cotensorbifunctor} is well
defined. In the special case of $\coring{C}=A$, we have
$-\cotensor{\coring{C}}-=-\tensor{A}-.$
\medskip

The conjunction of \cite[Proposition 2.7]{Zarouali:2004} and the
following result generalizes \cite[Theorem
2.4]{Brzezinski/Gomez:2003}.

\begin{theorem}\label{Comatrixcoring}
Let $X\in{}^{\coring{C}}\mathcal{M}%
^{\coring{D}}$ and
$\Lambda\in{}^{\coring{D}}\mathcal{M}^\coring{C}$. Assume that at
least one of the following conditions holds
\begin{enumerate}[(1)] \item
\begin{enumerate}[(a)] \item $_A\coring{C}$,
$_B\coring{D}$, $\coring{C}_A$ and $\coring{D}_B$ are flat, \item
$_{\coring{C}}X$ and $_{\coring{D}}\Lambda$ are coflat; or
\end{enumerate}
\item $A$ and $B$ are von Neumann regular rings; or
\item $_A\coring{C}$,
$_B\coring{D}$ are flat, and $\coring{C}$ and $\coring{D}$ are
coseparable corings.
\end{enumerate}
If there exist bicolinear maps
\[
\psi:\coring{C}\rightarrow X\square_{\coring{D}}\Lambda\;\textrm{and
}\;\omega:\Lambda\square_{\coring{C}}X\rightarrow\coring{D}
\]
in $^{\coring{C}}%
\mathcal{M}^\coring{C}$ and $^{\coring{D}}%
\mathcal{M}^\coring{D}$ respectively, such that the diagrams
\begin{equation}\label{unitcounit}
\xymatrix{\Lambda\ar[rr]^\simeq \ar[d]_\simeq &&
\Lambda\cotensor{\coring{C}}\coring{C}
\ar[d]^{\Lambda\cotensor{\coring{C}}\psi}
\\\coring{D}\cotensor{\coring{D}}\Lambda &&
\Lambda\cotensor{\coring{C}}X\cotensor{\coring{D}}\Lambda
\ar[ll]^{\omega\cotensor{\coring{D}}\Lambda}}
\xymatrix{X\ar[rr]^\simeq \ar[d]_\simeq &&
\coring{C}\cotensor{\coring{C}}X \ar[d]^{\psi\cotensor{\coring{C}}X}
\\X\cotensor{\coring{D}}\coring{D} &&
X\cotensor{\coring{D}}\Lambda\cotensor{\coring{C}}X
\ar[ll]^{X\cotensor{\coring{D}}\omega}}
\end{equation}
commute, then $\Lambda\square_{\coring{C}}X$ is a $B$-coring with
coproduct
\[
\xymatrix{\Delta:\Lambda\square_{\coring{C}}X \ar[r]^-\simeq &
\Lambda\square_{\coring{C}}(
\coring{C}\square_{\coring{C}}X)\ar[rr]^-{\Lambda\square
_{\coring{C}}(\mathcal{\psi}\square_{\coring{C}}X)} & &
\Lambda\square_{\coring{C}}((X\square_{\coring{D}}\Lambda)
\square_{\coring{C}}X)\ar[r]^-\simeq &
(\Lambda\square_{\coring{C}}X)\square_{\coring{D}}(
\Lambda\square_{\coring{C}}X) \ar@{^{(}->}[d] \\ & & & &
(\Lambda\square_{\coring{C}}X)\otimes_B(\Lambda
\square_{\coring{C}}X),}
 \]
 and counit
$\xymatrix{\epsilon=\epsilon_{\coring{D}}\circ\omega
:\Lambda\square_{\coring{C}}X\ar[r] & B}$.
\end{theorem}

\begin{proof}
By the bicolinearity of $\psi$, we have the commutativity of the two
diagrams
\begin{equation}\label{symetrydiagram1}
\xymatrix{\coring{C}\ar[r]^{\Delta_{\coring{C}}}
\ar[ddrr]^{(\psi\square_{\coring{C}}\psi)
\circ\Delta_{\coring{C}}} \ar[d]_{\psi} &
\coring{C}\square_{\coring{C}}\coring{C}
\ar[r]^-{\psi\square_{\coring{C}}\coring{C}}
& (X\square_{\coring{D}}\Lambda)\square_{\coring{C}}\coring{C}
\ar[dd]^{ (X\square_{\coring{D}}\Lambda)\square_{\coring{C}}\psi}
\\ (X\square_{\coring{D}}\Lambda)\ar[d]_\simeq & &
\\
(X\square_{\coring{D}}\Lambda)\square_{\coring{C}}\coring{C}
\ar[rr]_{(X\square_{\coring{D}}\Lambda) \square_{\coring{C}}\psi} &
& (X\square_{\coring{D}}\Lambda)
\square_{\coring{C}}(X\square_{\coring{D}}\Lambda)}
\end{equation}
\begin{equation}\label{symetrydiagram2}
\xymatrix{\coring{C}\ar[r]^{\Delta_{\coring{C}}}
\ar[ddrr]^{(\psi\square_{\coring{C}}\psi)
\circ\Delta_{\coring{C}}} \ar[d]_{\psi} &
\coring{C}\square_{\coring{C}}\coring{C}
\ar[r]^-{\coring{C}\square_{\coring{C}}\psi}
&
\coring{C}\square_{\coring{C}}(X\square_{\coring{D}}\Lambda)
\ar[dd]^{
\psi\square_{\coring{C}}(X\square_{\coring{D}}\Lambda)}
\\ (X\square_{\coring{D}}\Lambda)\ar[d]_\simeq & &
\\
\coring{C}\square_{\coring{C}}(X\square_{\coring{D}}\Lambda)
\ar[rr]_{\psi\square_{\coring{C}}
(X\square_{\coring{D}}\Lambda)} & &
(X\square_{\coring{D}}\Lambda)
\square_{\coring{C}}(X\square_{\coring{D}}\Lambda).}
\end{equation}
Therefore, we have the commutativity of the diagram
\[
\xymatrix{ \Lambda\square_{\coring{C}}X
\ar[r]^\simeq \ar[d]_\simeq&
\Lambda\square_{\coring{C}}\coring{C}
\square_{\coring{C}}X
\ar[rr]^{\Lambda\square_{\coring{C}}
\psi\square_{\coring{C}}X} & &
(\Lambda\square_{\coring{C}}X)\square_{\coring{D}}
(\Lambda\square_{\coring{C}}X)
\ar[d]^\simeq
 \\
\Lambda\square_{\coring{C}}\coring{C}
\square_{\coring{C}}X\ar[dd]_
{\Lambda\square_{\coring{C}}\psi
\square_{\coring{C}}X} & & &
(\Lambda\square_{\coring{C}}\coring{C}
\square_{\coring{C}}X)\square_{\coring{D}}
(\Lambda\square_{\coring{C}}X)
\ar[dd]|{(\Lambda\square_{\coring{C}}\psi
\square_{\coring{C}}X)\square_{\coring{D}}
(\Lambda\square_{\coring{C}}X)}
\\ & & & &
\\
(\Lambda\square_{\coring{C}}X)\square_{\coring{D}}
(\Lambda\square_{\coring{C}}X)
\ar[dr]_\simeq& & &
(\Lambda\square_{\coring{C}}X)\square_{\coring{D}}
(\Lambda\square_{\coring{C}}X)
\square_{\coring{D}}(\Lambda\square_{\coring{C}}X)
\\ &
(\Lambda\square_{\coring{C}}X)\square_{\coring{D}}
(\Lambda\square_{\coring{C}}\coring{C}\square_{\coring{C}}X)
\ar[urr]|{(\Lambda\square_{\coring{C}}X)\square_{\coring{D}}
(\Lambda\square_{\coring{C}}\psi\square_{\coring{C}}X)} & & &.}
\]
Hence the coassociative property of $\Delta$ follows. On the other
hand, if we put
$i:(\Lambda\square_{\coring{C}}X)\square_{\coring{D}}
(\Lambda\square_{\coring{C}}X)
\hookrightarrow(\Lambda\square_{\coring{C}}X)\otimes_B(\Lambda
\square_{\coring{C}}X)$ the canonical injection, we have,
\begin{multline}\label{6}
[( \epsilon_{\coring{D}}\circ
\omega\otimes_B(\Lambda\square_{\coring{C}}X)] \circ
i\circ(\Lambda\square_{\coring{C}}\mathcal{\psi}
\square_{\coring{C}}X)\circ[\xymatrix@1{\Lambda\square_{\coring{C}
}X\ar[r]^-\simeq &
\Lambda\square_{\coring{C}}\coring{C}\square_{\coring{C}}X}]
\\
 =[\epsilon_{\coring{D}}
\otimes_B(\Lambda\square_{\coring{C}}X)]\circ
i'\circ(\omega\square_{\coring{D}}(\Lambda
\square_{\coring{C}}X))\circ(\Lambda\square
_{\coring{C}}\mathcal{\psi}\square_{\coring{C}}X)\circ[\xymatrix@1{
\Lambda\square_{\coring{C}}X \ar[r]^-\simeq & \Lambda
\square_{\coring{C}}\coring{C}\square_{\coring{C}}X}],
\end{multline}
where $i':\coring{D}\square_{\coring{D}}(\Lambda\square
_{\coring{C}}X)
\hookrightarrow\coring{D}\otimes_B(\Lambda\square_{\coring{C}}X) $
is the canonical injection.

The first diagram of \eqref{unitcounit} is commutative means that
the diagram
\[
\xymatrix{\Lambda\square_{\coring{C}}
\coring{C}\ar[rr]^-{\Lambda\square
_{\coring{C}}\mathcal{\psi}} \ar[drr]_\simeq & & \Lambda\square
_{\coring{C}}(X\square_{\coring{D}}\Lambda)
\ar[rr]^-{\omega\square_{\coring{D}}\Lambda}
& & \coring{D}\square_{\coring{D}}\Lambda \\
 &  & \Lambda\ar[urr]_\simeq &  & }
\]
commutes. The composition of morphisms
\[
\xymatrix@1{\Lambda\ar[r]^-\simeq &
\Lambda\square_{\coring{C}}\coring{C}\ar[r]^-\simeq &
\Lambda\ar[r]^-\simeq &
\coring{D}\square_{\coring{D}}\Lambda\ar@{^{(}->}[r] &
\coring{D}\otimes_B\Lambda
\ar[rr]^{\epsilon_{\coring{D}}\otimes_B\Lambda} & &
B\otimes_B\Lambda}
\]
is exactly the morphism $[\lambda\mapsto1_B\otimes_B\lambda].$ Then,
\[
(\ref{6})=\Big[\sum_i\lambda_i\otimes_Ax_i\mapsto1_B\otimes
_B(\sum_i\lambda_i\otimes_Ax_i)\Big]  .
\]
Analogously, by using the commutativity of the second diagram of
\eqref{unitcounit}, we complete the proof of the counit property.
\end{proof}

\medskip

Let $N\in{}^{\coring{C}}\mathcal{M}^\coring{D}$ be a bicomodule,
quasi-finite as a right $\coring{D}$-comodule, such that
$_B\coring{D}$ is flat or $\coring{C}$ is a coseparable $A$-coring.
Let $T$ be a $k$-algebra.
\\ We will give an other characterization of the
natural isomorphism (for its definition we refer to
\cite{Zarouali:2004}):
\[
\Upsilon_{-,-}:-\otimes_{T}\cohom{\coring{D}}{N}{-}\rightarrow
\cohom{\coring{D}}{N}{-\otimes_{T}-}
\]
associated to the cohom functor $\cohom{\coring{D}}{N}{-}
:\mathcal{M}^\coring{D}\rightarrow\mathcal{M}^\coring{C}.$
\\ Let
$\theta:1_{\mathcal{M}%
^{\coring{D}}}\rightarrow h_{\coring{D}}(N,-)\otimes _AN$ be the
unit of the adjunction $(\cohom{\coring{D}}{N}{-},-\otimes_AN)$.
\\ Let $M\in{}^{T}\mathcal{M}^\coring{D}$ and
$W\in\mathcal{M}^{T}.$ Since the functors
$\cohom{\coring{D}}{N}{-}:\mathcal{M}%
^{\coring{D}}\rightarrow\mathcal{M}^{A}$ and
$-\otimes_AN:\mathcal{M}%
^{A}\rightarrow\mathcal{M}^\coring{D}$ are $k$-linear and preserve
inductive limits, then the functor $\cohom{\coring{D}}{N}{-}
\otimes_AN:\mathcal{M}^\coring{D}\rightarrow\mathcal{M}^{\coring{D}%
}$ is also $k$-linear and preserves inductive limits. By \cite[Lemma
3.2(1)]{Gomez:2002}, $\theta_M$ is a morphism in
$^{T}\mathcal{M}^\coring{D}.$ Define
\[
\Xi_{W,M}:\cohom{\coring{D}}{N}{W\otimes_{T}M}\rightarrow
W\otimes_{T}\cohom{\coring{D}}{N}{M}
\]
to be the unique morphism in $\mathcal{M}^{A}$ satisfying $(\Xi
_{W,M}\otimes_AN) \theta_{W\otimes_{T}M}=W\otimes_T\theta_M$. We
have
\begin{multline}
\Big[\rho_{W\otimes_T\cohom{\coring{D}}{N}{M}}\otimes_AN-\big(W
\otimes_T\cohom{\coring{D}}{N}{M}\big)
\otimes_A\lambda_N\Big](W\otimes_{T}\theta_M)   \\
=W\otimes_{T}\Big[\big(\rho_{\cohom{\coring{D}}{N}{M}
}\otimes_AN-\cohom{\coring{D}}{N}{M}\otimes_A\lambda_N\big)
\theta_M\Big]=0.
\end{multline}
Then
$\operatorname{Im}(W\otimes_T\theta_M)\subset(W\otimes_T\cohom{\coring{D}}{N}{M}
\square_{\coring{C}}N,$ and by the proof of \cite[Proposition
4.2(1)]{Gomez:2002}, $\Xi_{W,M}$ is a morphism in
$\mathcal{M}^\coring{C}$. Now we will verify that $\Xi_{-,M}$ is a
natural transformation. For this, let $f:W\rightarrow W'$ be a
morphism in $\mathcal{M}^{T}.$ We will verify the commutativity of
the diagram
\[%
\xymatrix{\cohom{\coring{D}}{N}{W\otimes_TM}\ar[d]
\ar[rr]^{\Xi_{W,M}
} & & W\otimes_T\cohom{\coring{D}}{N}{M} \ar[d] \\
 \cohom{\coring{D}}{N}{W'\otimes_TM}
\ar[rr]^{\Xi_{W',M}} & &
 W'\otimes_T\cohom{\coring{D}}{N}{M}.}
\]
Since $\theta$ is a natural transformation, the diagram
\[%
 \xymatrix{W\otimes_TM \ar[d]
\ar[rr]^-{\theta_{W\otimes_TM}} & &
\cohom{\coring{D}}{N}{W\otimes_TM}\otimes_AN \ar[d]\\
W'\otimes_TM \ar[rr]^-{\theta_{W'\otimes_TM}} & &
\cohom{\coring{D}}{N}{W'\otimes_TM}\otimes_AN}
\]
commutes. By a straightforward computations we have
$$\Big[\big(\Xi_{W',M}\cohom{\coring{D}}{N}{f\otimes_TM}
\big)\otimes_AN\Big]\theta_{W\otimes_TM}\newline =\bigg[\Big(\big(
f\otimes_T\cohom{\coring{D}}{N}{M}\big)\Xi_{W,M}\Big)\otimes_AN\bigg]
\theta_{W\otimes_TM}$$ By uniqueness,
$\Xi_{W',M}\cohom{\coring{D}}{N}{f\otimes_TM}=(
f\otimes_T\cohom{\coring{D}}{N}{M})\Xi_{W,M}$. Finally, we will
verify that $\Upsilon_{T,M}=\Xi_{T,M}^{-1}$, i.e. $\Xi_{T,M}$
satisfies the commutativity of the diagram
\[%
\xymatrix{\cohom{\coring{D}}{N}{T\otimes_TM}\ar[rd]_\simeq
\ar[rr]^{\Xi_{T,M} }&  &
T\otimes_T\cohom{\coring{D}}{N}{M}\ar[ld]^\simeq
\\
 &  \cohom{\coring{D}}{N}{M} & }
\]
(the canonical isomorphism $\xymatrix{T\otimes_TM \ar[r]^-\simeq &
M}$ is an isomorphism in $\mathcal{M}^\coring{D}$). Since $\theta$
is a natural transformation, the diagram
\[%
\xymatrix{T\otimes_TM \ar[d]^\simeq \ar[rr]^-{\theta_{T\otimes_TM}}
& & \cohom{\coring{D}}{N}{T\otimes_TM}\otimes_AN
\ar[d]^\simeq \\
 M \ar[rr]^-{\theta_M} &  &
\cohom{\coring{D}}{N}{M}\otimes_AN}
\]
commutes. It follows from the commutativity of the last diagram and
the fact that $\theta_M$ is $T$-linear, that the following diagram
is also commutative
\[%
\xymatrix{ T\otimes_TM \ar[d]^\simeq \ar[rr]^-{\theta_{T\otimes_TM}}
& & \cohom{\coring{D}}{N}{T\otimes_TM} \otimes_AN
\ar[d]^\simeq\\
 T\otimes_T\cohom{\coring{D}}{N}{M}\otimes_AN
\ar[rr]^\simeq  & & \cohom{\coring{D}}{N}{M}\otimes_AN.}
\]
By \cite[Corollary 3.6.6]{Popescu:1973}, $\Xi_{W,M}$ is an
isomorphism, and by a theorem of Mitchell (\cite[Theorem
3.6.5]{Popescu:1973}), $\Upsilon_{W,M}=\Xi_{W,M}^{-1}$.

\medskip
We obtain then the following generalization of
\cite[1.6]{Takeuchi:1977}:

\begin{proposition}\label{Upsilon}
Let $N\in{}^{\coring{C}}\mathcal{M}^\coring{D}$ be a bicomodule,
quasi-finite as a right $\coring{D}$-comodule, such that
$_B\coring{D}$ is flat or $\coring{C}$ is a coseparable $A$-coring.
Let
$\theta:1_{\mathcal{M}%
^{\coring{D}}}\rightarrow \cohom{\coring{D}}{N}{-}\otimes_AN$ be the
unit of the adjunction $(\cohom{\coring{D}}{N}{-},-\otimes_AN)$, and
let $T$ be a $k$-algebra. If $\Upsilon
_{-,-}:-\otimes_T\cohom{\coring{D}}{N}{-}\rightarrow
\cohom{\coring{D}}{N}{-\otimes_T-}$ is the natural isomorphism
associated to the cohom functor
$\cohom{\coring{D}}{N}{-}:\mathcal{M}^\coring{D}\to\mathcal{M}^\coring{C},$
then $\Upsilon_{W,M}=\Xi_{W,M}^{-1},$ where $\Xi_{W,M}$ is defined
as above.
\end{proposition}

\medskip
As a consequence of the last result we get the following
generalization of \cite[1.13]{Takeuchi:1977}.

\begin{proposition}\label{delta}
Let $N\in\bcomod{\coring{C}}{\coring{D}}$ be a bicomodule,
quasi-finite as a right $\coring{D}$-comodule, such that
$_B\coring{D}$ is flat. Suppose that $_A\coring{C}$ is flat (resp.
 $\mathcal{M}^\coring{C}$ is an abelian category). Let
$\theta:1_{\mathcal{M}%
^{\coring{D}}}\rightarrow \cohom{\coring{D}}{N}{-}\otimes_AN$ be the
unit of the adjunction $(\cohom{\coring{D}}{N}{-},-\otimes_AN)$. If
the cohom functor $\cohom{\coring{D}}{N}{-}$ is exact (resp.
$\coring{D}$ is a coseparable $B$-coring), then the natural
isomorphism $\delta$ (see \cite[Corollary 3.5]{Zarouali:2004}):
\[
\cohom{\coring{D}}{N}{-}
\simeq-\square_{\coring{D}}\cohom{\coring{D}}{N}{\coring{D}}:
\mathcal{M}^\coring{D}\rightarrow \mathcal{M}^\coring{C}
\]
satisfy the property:\\
For every $M\in\mathcal{M}^\coring{D}$,
$\delta_M:\cohom{\coring{D}}{N}{M}\rightarrow
M\square_{\coring{D}}\cohom{\coring{D}}{N}{\coring{D}}$ is the
unique right $A$-linear map satisfying
\[
(\delta_M\otimes_AN)\theta_M=(M\square_{\coring{D}}
\theta_{\coring{D}})\rho_M.
\]
\end{proposition}

\begin{proof}
 Let $M\in\mathcal{M}^\coring{D}$, we have
 $(\Xi_{M,\coring{D}}\otimes_AN)\theta_{M\otimes_B\coring{D}}=
 M\otimes_B\theta_{\coring{D}}$. On the other hand,
since $\theta$ is a natural transformation,
$(\Xi_{M,\coring{D}}\otimes_AN)
 \theta_{M\otimes_B\coring{D}}\rho_M=(\Xi_{M,\coring{D}}\otimes_AN)
 (\cohom{\coring{D}}{N}{\rho_M}\otimes_AN)
 \theta_M=(i\delta_M\otimes_AN)\theta_M$, where
 $i:M\square_{\coring{D}}\cohom{\coring{D}}{N}{\coring{D}}
 \hookrightarrow
M\otimes_B\cohom{\coring{D}}{N}{\coring{D}}$ is the canonical
injection. Therefore,
$(i\otimes_AN)(\delta_M\otimes_AN)\theta_M=(M\otimes_B
\theta_{\coring{D}})\rho_M$.
 Since $i\otimes_AN$ is a monomorphism ($_AN$ is
flat), $(\delta_M\otimes_AN)\theta_M
 =(M\square_{\coring{D}}\theta_{\coring{D}})\rho_M$.
\end{proof}

\medskip

Now we will give a series of properties concerning our comatrix
coring.

\begin{proposition}\label{comatrixproperty1}
Let $\Lambda\in \bcomod{\coring{D}}{\coring{C}}$ be a bicomodule,
quasi-finite as a right $\coring{C}$-comodule, such that
$_A\coring{C}$ and $_B\coring{D}$ are flat. Set
$X=\cohom{\coring{C}}{\Lambda}{\coring{C}}\in
\bcomod{\coring{D}}{\coring{C}}$.
 If
 \begin{enumerate}[(a)]\item $\coring{C}_A$ and $\coring{D}_B$ are flat,
 the cohom functor $\cohom{\coring{C}}{\Lambda}{-}$ is exact and
$_{\coring{D}}{\Lambda}$ is coflat, or
\item $A$ and $B$ are von Neumann regular rings and the cohom functor
$\cohom{\coring{C}}{\Lambda}{-}$ is exact, or
\item $\coring{C}$ and
$\coring{D}$ are coseparable corings,
\end{enumerate}
then we have
\begin{enumerate} [(1)]\item
$\xymatrix{\delta_{\Lambda}:\e{\coring{C}}{\Lambda}\ar[r] &
\Lambda\cotensor{\coring{C}}X}$ is an isomorphism of $B$-corings;
\item
$\xymatrix{\chi_{\coring{D}}\circ
\cohom{\coring{C}}{\Lambda}{\lambda_{\Lambda}}=
\omega\circ\delta_{\Lambda}: \e{\coring{C}}{\Lambda}\ar[r] &
\coring{D}}$ and
$\xymatrix{\omega:\Lambda\cotensor{\coring{C}}X\ar[r] &\coring{D}}$
is a homomorphism of $B$-corings, where
$\lambda_{\Lambda}:\Lambda\to\coring{D}\cotensor{\coring{D}}{\Lambda}$
is the left coaction of $\Lambda$, $\chi$ the counit of the
adjunction
$(\cohom{\coring{C}}{\Lambda}{-},-\square_{\coring{D}}\Lambda)$, and
$\omega$ is the unique $\coring{D}-\coring{D}$-bicolinear map such
that \eqref{counitcohom} is the counit of the adjunction
$(-\cotensor{\coring{C}}{\cohom{\coring{C}}{\Lambda}{\coring{C}}},
-\cotensor{\coring{D}} {\Lambda})$.
\end{enumerate}
\end{proposition}

\begin{proof}
Let
$\theta:1_{\mathcal{M}%
^{\coring{C}}}\rightarrow
\cohom{\coring{C}}{\Lambda}{-}\tensor{B}{\Lambda}$ be the unit of
the adjunction
$(\cohom{\coring{C}}{\Lambda}{-},-\tensor{B}{\Lambda})$.

(1) By \cite[Theorem 2.3]{Zarouali:2004},
$\cohom{\coring{C}}{\Lambda}{-}\simeq
-\square_{\coring{C}}\cohom{\coring{C}}{\Lambda}{\coring{C}}$. In
the first case, $\cohom{\coring{C}}{\Lambda}{\coring{C}}$ is coflat
as a left $\coring{C}$-comodule. By \cite[Proposition
3.4]{Gomez:2002}, $\theta_{\coring{C}}:\coring{C}\rightarrow
\cohom{\coring{C}}{\Lambda}{\coring{C}} \square_{\coring{D}}\Lambda$
is a $\coring{C}$-bicolinear
map.\\
From Proposition \ref{delta}, for every
$M\in{}\mathcal{M}^\coring{C}$, we have a commutative diagram
\[
\xymatrix{M\ar[d]_{\theta_M}\ar[rr]^\simeq & &
M\square_{\coring{C}}\coring{C}\ar[d]^{M
\square_{\coring{C}}\theta_{\coring{C}}}
\\
\cohom{\coring{C}}{\Lambda}{M}\square_{\coring{D}}\Lambda
\ar[rr]_-{\delta_M \square_{\coring{D}}\Lambda} & &
M\square_{\coring{C}}\cohom{\coring{C}}{\Lambda}
{\coring{C}}\square_{\coring{D}}\Lambda.}
\]
From \cite[Lemma 2.6]{Zarouali:2004}, the unit and the counit of the
adjunction
$(-\cotensor{\coring{C}}{\cohom{\coring{C}}
{\Lambda}{\coring{C}}},-\cotensor{\coring{D}}
{\Lambda})$ are
\begin{equation}\label{unitcohom}
\xymatrix@1{\eta:1_{\mathcal{M}^\coring{C}}\ar[r]^\simeq &
-\square_{\coring{C}}\coring{C}\ar[rr]^-
{-\square_{\coring{C}}\theta_{\coring{C}}}
& &
-\square_{\coring{C}}\cohom{\coring{C}}{\Lambda}{\coring{C}}
\square_{\coring{D}}\Lambda}
\end{equation}
 and
\begin{equation}\label{counitcohom}
\xymatrix@1{\varepsilon:-\square_{\coring{D}}
\Lambda\square_{\coring{C}}
\cohom{\coring{C}}{\Lambda}{\coring{C}}
\ar[rr]^-{-\square_{\coring{D}}\omega} & &
-\square_{\coring{D}}\coring{D}\ar[r]^\simeq &
1_{\mathcal{M}^\coring{D}}.}
\end{equation}
To show that $\delta_\Lambda$ is a homomorphism of $B$-corings, it
suffices to prove (see \cite[23.8]{Brzezinski/Wisbauer:2003}):
\begin{multline}
(\delta_{\Lambda}\tensor{B}{\delta_{\Lambda}}\tensor{B}{\Lambda})
\circ(\e{\coring{C}}{\Lambda}\tensor{B}{\theta_{\Lambda}})
\circ\theta_{\Lambda}
\\
=(\Lambda\cotensor{\coring{C}}{\theta_\coring{C}}
\cotensor{\coring{C}}X\tensor{B}{\Lambda})
\circ[\xymatrix{\Lambda\cotensor{\coring{C}}X\tensor{B}{\Lambda}
\ar[r]^-\simeq&
\Lambda\cotensor{\coring{C}}{\coring{C}}\cotensor{\coring{C}}X
\tensor{B}{\Lambda}}]
\circ(\delta_{\Lambda}\tensor{B}{\Lambda})\circ\theta_{\Lambda},
\end{multline}
and
\begin{multline}
(\epsilon_{\coring{D}}\tensor{B}{\Lambda})\circ
(\omega\tensor{B}{\Lambda})
\circ(\delta_{\Lambda}\tensor{B}{\Lambda})
\circ\theta_{\Lambda}=
[\xymatrix{\Lambda\ar[r]^-\simeq& B\tensor{B}{\Lambda}}].
\end{multline}

Now we will prove $(8)$. We have,
\begin{multline}
\begin{split}
(\delta_{\Lambda}\tensor{B}{\delta_{\Lambda}}
\tensor{B}{\Lambda})
\circ(\e{\coring{C}}{\Lambda}\tensor{B}
{\theta_{\Lambda}})
\circ\theta_{\Lambda}
&
=\Big[\delta_{\Lambda}\tensor{B}{\big((\delta_{\Lambda}
\tensor{B}{\Lambda})\circ\theta_{\Lambda}\big)}\Big]
\circ\theta_{\Lambda}\\
&
=\Big[(\Lambda\cotensor{\coring{C}}X)\tensor{B}
{\big((\delta_{\Lambda}\tensor{B}{\Lambda})
\circ\theta_{\Lambda}\big)\Big]}\circ(\delta_{\Lambda}
\tensor{B}{\Lambda})\circ\theta_{\Lambda}\\
& =\Big[(\Lambda\cotensor{\coring{C}}X)\tensor{B}
{\big((\Lambda\cotensor{\coring{C}}
{\theta_{\coring{C}}})\circ\rho_{\Lambda})\Big]}
\circ(\Lambda\cotensor{\coring{C}}
{\theta_{\coring{C}}})\circ\rho_{\Lambda}\\ &
\qquad\textrm{(using
Proposition \ref{delta})}\\
&
=(\Lambda\cotensor{\coring{C}}X\tensor{B}
{\Lambda\cotensor{\coring{C}}
{\theta_{\coring{C}}}})\circ
(\Lambda\cotensor{\coring{C}}X\tensor{B}
{\rho_{\Lambda}})\circ(\Lambda\cotensor{\coring{C}}
{\theta_{\coring{C}}})\circ\rho_{\Lambda}
\end{split}
\end{multline}
On the other hand, by Proposition \ref{delta},
\begin{multline}
(\Lambda\cotensor{\coring{C}}{\theta_\coring{C}}
\cotensor{\coring{C}}X\tensor{B}{\Lambda})
\circ[\xymatrix{\Lambda\cotensor{\coring{C}}X
\tensor{B}{\Lambda}
\ar[r]^-\simeq&
\Lambda\cotensor{\coring{C}}{\coring{C}}
\cotensor{\coring{C}}X\tensor{B}{\Lambda}}]
\circ(\delta_{\Lambda}\tensor{B}{\Lambda})
\circ\theta_{\Lambda}\\
=(\Lambda\cotensor{\coring{C}}{\theta_\coring{C}}
\cotensor{\coring{C}}X\tensor{B}{\Lambda})
\circ[\xymatrix{\Lambda\cotensor{\coring{C}}X\tensor{B}
{\Lambda}\ar[r]^-\simeq & \Lambda\cotensor{\coring{C}}
{\coring{C}}\cotensor{\coring{C}}X\tensor{B}{\Lambda}}]
\circ(\Lambda\cotensor{\coring{C}}{\theta_{\coring{C}}})
\circ\rho_{\Lambda}.
\end{multline}
 Now, by \cite[Proposition 3.4]{Gomez:2002},
$\theta_{\coring{C}}:\coring{C}\rightarrow
\cohom{\coring{C}}{\Lambda}{\coring{C}}\tensor{B}{\Lambda}$ is
 $\coring{C}$-bicolinear. Then
$\theta_{\coring{C}}$
 makes commutative the diagrams
\eqref{symetrydiagram1} and \eqref{symetrydiagram2}
 (by replacing $\psi$ by $\theta_{\coring{C}}$). Hence
$(10)=(11)$. \\
 Finally, we will prove $(9)$.
\begin{multline*}
\begin{split}
(\epsilon_{\coring{D}}\tensor{B}{\Lambda})
\circ(\omega\tensor{B}{\Lambda})
\circ(\delta_{\Lambda}\tensor{B}{\Lambda})
\circ\theta_{\Lambda} &
=(\epsilon_{\coring{D}}\tensor{B}{\Lambda})
\circ(\omega\tensor{B}{\Lambda})
\circ(\Lambda\cotensor{\coring{C}}
{\theta_{\coring{C}}})\circ\rho_{\Lambda}
\\ & \qquad \textrm{(using Proposition \ref{delta})}\\
& =(\epsilon_{\coring{D}}\tensor{B}{\Lambda})
\circ\lambda_{\Lambda}
\qquad\textrm{(using the first equality of
\eqref{unitcounit})}\\
& =[\xymatrix{\Lambda \ar[r]^-\simeq&
B\tensor{B}{\Lambda}}].
\end{split}
\end{multline*}

(2) Since $\delta$ is a natural isomorphism, the diagram
\[
\xymatrix{\cohom{\coring{C}}{\Lambda}{\Lambda}
\ar[d]_{\cohom{\coring{C}}{\Lambda}{\lambda_{\Lambda}}}
\ar[rr]^{\delta_{\Lambda}}
& & {\Lambda}\cotensor{\coring{C}}X
\ar[d]^{\lambda_{\Lambda}\cotensor{\coring{C}}X}
\\
\cohom{\coring{C}}{\Lambda}{\coring{D}\cotensor{\coring{D}}{\Lambda}}
\ar[rr]_-{\delta_{\coring{D}\cotensor{\coring{D}}{\Lambda}}} & &
(\coring{D}\cotensor{\coring{D}}{\Lambda})\cotensor{\coring{C}}X}
\]
is commutative. The counit of the adjunction
$(\cohom{\coring{C}}{\Lambda}{-},-\cotensor{\coring{D}}{\Lambda})$
is \begin{equation}
\xymatrix@1{\cohom{\coring{C}}{\Lambda}{-\cotensor{\coring{D}}
{\Lambda}}\ar[rr]^{\delta_{-\cotensor{\coring{D}}{\Lambda}}} & &
-\cotensor{\coring{D}}{\Lambda}
\cotensor{\coring{C}}X\ar[rr]^-{-\cotensor{\coring{D}}{\omega}} & &
-\cotensor{\coring{D}}{\coring{D}}\ar[r]^\simeq &
1_{\mathcal{M}^\coring{D}}.}
\end{equation}
Then,
$\chi_{\coring{D}}\circ\cohom{\coring{C}}{\Lambda}{\lambda_{\Lambda}}=
[\xymatrix{\coring{D}\cotensor{\coring{D}}{\coring{D}}\ar[r]^\simeq
& \coring{D}}]\circ(\coring{D}\cotensor{\coring{D}}{\omega})
\circ(\lambda_{\Lambda}\cotensor{\coring{C}}X)\circ\delta_{\Lambda}$.\\
Since $\omega$ is left $\coring{D}$-colinear, and
$\lambda_{\Lambda\cotensor{\coring{C}}X}=\lambda_{\Lambda}
\cotensor{\coring{C}}X$,
and from the commutativity of the following diagram
\[
\xymatrix{\Lambda\cotensor{\coring{C}}X \ar[r]^-\simeq
\ar[dr]_{\lambda_{\Lambda\cotensor{\coring{C}}X}} &
\coring{D}\cotensor{\coring{D}}(\Lambda\cotensor{\coring{C}}X)
\ar[rr]^{\coring{D}\cotensor{\coring{D}}{\omega}} \ar@{^{(}->}[d] &
&
\coring{D}\cotensor{\coring{D}}\coring{D}\ar[r]^-\simeq\ar@{^{(}->}[d]
 & \coring{D} \ar[dl]^{\Delta_{\coring{D}}}
 \\ &
\coring{D}\tensor{B}(\Lambda\cotensor{\coring{C}}X)
\ar[rr]_{\coring{D}\tensor{B}{\omega}}
 & & \coring{D}\tensor{B}{\coring{D}} & ,}
\]
$\chi_{\coring{D}}\circ\cohom{\coring{C}}{\Lambda}{\lambda_{\Lambda}}=
\omega\circ\delta_{\Lambda}$. Finally, by \cite[Proposition
5.2]{ElKaoutit/Gomez:2003}, it is a homomorphism of $B$-corings.
Hence, $\omega$ is also a homomorphism of $B$-corings.
\end{proof}

\begin{example}
Let $_A\coring{C}$ and $_B\coring{D}$ be flat. Let
$(\varphi,\rho):\coring{C}\rightarrow\coring{D}$ be a homomorphism
of corings. Suppose that $\coring{C}_A$ and $\coring{D}_B$ are flat
and $_{\coring{D}}(B\tensor{A}{\coring{C}})$ is coflat, or
$\coring{C}$ and $\coring{D}$ are coseparable corings. Then, we have
$(-\tensor{A}{B},-\cotensor{\coring{D}}(B\tensor{A}{\coring{C}}))$
is an adjoint pair, and
$-\tensor{A}{B}\simeq-\cotensor{\coring{C}}{(\coring{C}\tensor{A}{B})}$
(see the proof of \cite[Theorem 4.1]{Zarouali:2004}). By
\cite[23.9]{Brzezinski/Wisbauer:2003}, the comatrix coring
$(B\tensor{A}{\coring{C}})\cotensor{\coring{C}}(\coring{C}\tensor{A}{B})$
($\simeq\e{\coring{C}}{B\tensor{A}{\coring{C}}}$) is isomorphic (as
corings) to the coring $B\coring{C}B$ (see
\cite[17.2]{Brzezinski/Wisbauer:2003}).
\end{example}

\begin{theorem}\label{comatrixproperty2}
Let $_A\coring{C}$ be flat. Let $\Lambda
\in{}_B\mathcal{M}^\coring{C}$ be a bicomodule, quasi-finite as a
right $\coring{C}$-comodule. Set
$X=\cohom{\coring{C}}{\Lambda}{\coring{C}}\in{}^{\coring{C}}\mathcal{M}_B$.
 If $\coring{C}_A$ is flat and the cohom functor
$\cohom{\coring{C}}{\Lambda}{-}$ is exact, or if $\coring{C}$ is a
coseparable coring, then we have
\begin{enumerate} [(1)]\item The functor
$-\tensor{B}{\Lambda}:
\mathcal{M}_B\rightarrow\mathcal{M}^\coring{C}$ is separable if and
only if the comatrix coring $\Lambda\cotensor{\coring{C}}X$ is a
cosplit $B$-coring. \item If the cohom functor
$\cohom{\coring{C}}{\Lambda}{-}$ is separable then the comatrix
coring $\Lambda\cotensor{\coring{C}}X$ is a coseparable $B$-coring.
\item If $-\tensor{B}{\Lambda}:
\mathcal{M}_B\rightarrow\mathcal{M}^\coring{C}$ is a Frobenius
functor, either $X_B$ is flat and $\Lambda_{\coring{C}}$ is coflat
or $\coring{C}$ is coseparable, and if for all
$f\in(\Lambda\cotensor{\coring{C}}X)^*$ and all $b\in B$,
$bH(f)=H(f)b$, where $\xymatrix{H:(\Lambda\cotensor{\coring{C}}X)^*
\ar[r] & {\Lambda\cotensor{\coring{C}}X}}$ is the isomorphism
defined by \eqref{frobeniuscomatrix} (this is the case if
$\Lambda\cotensor{\coring{C}}X$ is a coalgebra (i.e. $B=k$)), then
the comatrix coring $\Lambda\cotensor{\coring{C}}X$ is a Frobenius
$B$-coring.
\end{enumerate}
\end{theorem}

\begin{proof}
(1) Let $\omega:\Lambda\cotensor{\coring{C}}X\rightarrow B$ the
unique $B$-bilinear map such that
\begin{equation}
\xymatrix@1{\varepsilon:-\otimes_B\Lambda\square_{\coring{C}}
\cohom{\coring{C}}{\Lambda}{\coring{C}} \ar[rr]^-{-\otimes_B\omega}
& & -\otimes_B\coring{D}\ar[r]^\simeq & 1_{\mathcal{M}_B}}
\end{equation}
is the counit of the adjunction
$(-\cotensor{\coring{C}}{\cohom{\coring{C}}{\Lambda}
{\coring{C}}},-\tensor{B}{\Lambda})$. We have,
$\epsilon_{\Lambda\cotensor{\coring{C}}X}=\omega$. By \cite[Lemma
2.6]{Zarouali:2004} and Rafael's theorem \cite [Theorem
24]{Caenepeel/Militaru/Zhu:2002}, the functor $-\tensor{B}{\Lambda}:
\mathcal{M}_B\rightarrow\mathcal{M}^\coring{C}$ is separable if and
only if there exists a $B$-bilinear map
$\omega':B\rightarrow\Lambda\cotensor{\coring{C}}X$
 satisfying $\omega\circ\omega'=1_B$. Hence, (1)
 follows from the definition of a cosplit coring (see
\cite{Brzezinski/Wisbauer:2003}).

(2) We know that the unit of the adjunction
$(-\cotensor{\coring{C}}{\cohom{\coring{C}}{\Lambda}
{\coring{C}}},-\tensor{B} {\Lambda})$ is
\begin{equation}
\xymatrix@1{\eta:1_{\mathcal{M}^\coring{C}}\ar[r]^\simeq &
-\square_{\coring{C}}\coring{C}\ar[rr]^-{-\square_{\coring{C}}
\theta_{\coring{C}}}
& &
-\square_{\coring{C}}\cohom{\coring{C}}{\Lambda}
{\coring{C}}\otimes_B\Lambda}.
\end{equation}
By \cite[Lemma 2.6]{Zarouali:2004} and Rafael's theorem \cite
[Theorem 24]{Caenepeel/Militaru/Zhu:2002}, the cohom functor is
separable if and only if there exists a
$\coring{C}-\coring{C}$-bicolinear map
$\psi':\cohom{\coring{C}}{\Lambda}{\coring{C}}\otimes_B\Lambda\to
\coring{C}$ satisfying
$\psi'\circ\theta_{\coring{C}}=1_{\coring{C}}$. We have,
$\Delta_{\Lambda\cotensor{\coring{C}}X}=(\Lambda
\cotensor{\coring{C}}{\theta_{\coring{C}}}
\cotensor{\coring{C}}X)\circ
[\xymatrix{\Lambda\cotensor{\coring{C}}X \ar[r]^-\simeq &
\Lambda\cotensor{\coring{C}}{\coring{C}}\cotensor
{\coring{C}}X}]$.\\
Let
$f=[\xymatrix{\Lambda\cotensor{\coring{C}}{\coring{C}}
\cotensor{\coring{C}}X\ar[r]^-\simeq
& \Lambda\cotensor{\coring{C}}X}]\circ
(\Lambda\cotensor{\coring{C}}{\psi'}\cotensor{\coring{C}}X)$. We
have,
$f\circ\Delta_{\Lambda\cotensor{\coring{C}}X}=1_{\Lambda
\cotensor{\coring{C}}X}$,
\begin{multline*}
\begin{split}
 \Delta_{\Lambda\cotensor{\coring{C}}X}\circ f &
=\big[[\xymatrix{\Lambda\cotensor{\coring{C}}
{\coring{C}}\ar[r]^-\simeq
&
\Lambda}]\circ(\Lambda\cotensor{\coring{C}}{\psi'})
\big]\cotensor{\coring{C}}
\big[(\psi\cotensor{\coring{C}}X)
\circ[\xymatrix{X\ar[r]^-\simeq &
\coring{C}\cotensor{\coring{C}}X}]\big]\\ &
=(f\tensor{B}{\Lambda\cotensor{\coring{C}}X})\circ
(\Lambda\cotensor{\coring{C}}X\tensor{B}
{\Delta_{\Lambda\cotensor{\coring{C}}X}}),
\end{split}
\end{multline*}
and
\begin{multline*}
\begin{split}
 \Delta_{\Lambda\cotensor{\coring{C}}X}\circ f &
=\big[(\Lambda\cotensor{\coring{C}}{\psi})
\circ[\xymatrix{\Lambda\ar[r]^-\simeq
 & \Lambda\cotensor{\coring{C}}{\coring{C}}}]
 \big]\cotensor{\coring{C}}
\big[[\xymatrix{\coring{C}\cotensor{\coring{C}}X
\ar[r]^-\simeq &
 X}]\circ(\psi'\cotensor{\coring{C}}X)\big] \\
 & =(\Lambda\cotensor{\coring{C}}X\tensor{B}{f})\circ
(\Delta_{\Lambda\cotensor{\coring{C}}X}\tensor{B}
{\Lambda\cotensor{\coring{C}}X}).
\end{split}
\end{multline*}
From the definition of a coseparable coring (see \cite{Guzman:1989}
or \cite{Brzezinski/Wisbauer:2003}), $\Lambda\cotensor{\coring{C}}X$
is a coseparable coring.

(3) We will prove it using \cite[27.13]{Brzezinski/Wisbauer:2003}.
Set $T=(\Lambda\cotensor{\coring{C}}X)^*$,
$\Delta=\Delta_{\Lambda\cotensor{\coring{C}}X}$ and
$\epsilon=\epsilon_{\Lambda\cotensor{\coring{C}}X}$. At first from
\cite[Theorem 2.11]{Zarouali:2004},
$(X\tensor{B}{-},\Lambda\cotensor{\coring{C}}{-})$ is a Frobenius
pair, and then $\Lambda\cotensor{\coring{C}}X\tensor{B}{-}$ is a
Frobenius functor. Hence $(\Lambda\cotensor{\coring{C}}X)_B$ is
finitely generated projective.\\
Let us consider the isomorphism
$$\xymatrix{H:T \ar[r]^-\simeq &
{\e{\coring{C}}{\Lambda}^*} \ar[r]^-{\phi_{\Lambda,B}}&
{\operatorname{End}_{\coring{C}}(\Lambda)} \ar[r]^-\simeq &
{\hom{\coring{C}}{B\tensor{B}{\Lambda}}{\Lambda}} \ar[r]^-\simeq &
{\hom{B}{B}{\Lambda\cotensor{\coring{C}}X}} \ar[r]^-\simeq &
{\Lambda\cotensor{\coring{C}}X}},$$ where
$\xymatrix{\phi_{\Lambda,B}:\e{\coring{C}}{\Lambda}^*\ar[r] &
{\operatorname{End}_{\coring{C}}(\Lambda)}}$ is the canonical
anti-isomorphism of rings defined in
\cite[23.8]{Brzezinski/Wisbauer:2003}. Then, (by using the following
consequence of Proposition \ref{delta};
$\Delta=\big[(\delta_{\Lambda}\tensor{B}{\Lambda})
\theta_{\Lambda}\big]\cotensor{\coring{C}}X$),
for $f\in T$,
\begin{multline}\label{frobeniuscomatrix}
H(f)=\Big[\big(\phi_{\Lambda,B}(f\delta_{\Lambda})
[\xymatrix{B\tensor{B}{\Lambda}\ar[r]^-\simeq
& {\Lambda}}]\big)\cotensor{\coring{C}}X\Big]\xi_B(1_B) \\
={[\xymatrix{B\tensor{B}{\Lambda}
\cotensor{\coring{C}}X\ar[r]^-\simeq
&
{\Lambda\cotensor{\coring{C}}X}}]
(f\tensor{B}{\Lambda}\cotensor{\coring{C}}X)
\Delta[\xymatrix{B\tensor{B}{\Lambda}
\cotensor{\coring{C}}X\ar[r]^-\simeq
& {\Lambda\cotensor{\coring{C}}X}}]}\xi_B(1_B),
\end{multline}
where
$\xi:1_{\mathcal{M}_B}\rightarrow-\tensor{B}
{\Lambda}\cotensor{\coring{C}}X$
is the unit of the adjunction
$(-\tensor{B}{\Lambda},-\cotensor{\coring{C}}X)$.

Now, let $f,f'\in T$,
$\phi_{\Lambda,B}(f*^{r}f')=\phi_{\Lambda,B}
(f'\delta_{\Lambda})
\phi_{\Lambda,B}(f'\delta_{\Lambda})$.
\begin{multline*}
\begin{split}
H(f*^{r}f') &
=\big(\phi_{\Lambda,B}(f'\delta_{\Lambda})
\cotensor{\coring{C}}X\big)\big(H(f)\big)
\\
&
=[\xymatrix{B\tensor{B}{\Lambda}\cotensor{\coring{C}}
X\ar[r]^-\simeq
&
{\Lambda\cotensor{\coring{C}}X}}](f'\tensor{B}{\Lambda}
\cotensor{\coring{C}}X)
\Delta\big(H(f)\big)\\
 \\
& =H(f).f'.
\end{split}
\end{multline*}
Finally, let $f\in T$, $b\in B$ and let $i_R:B\rightarrow T$ be the
anti-morphism of rings defined by $i_R(b)=\epsilon(b-)$ for $B\in B$
(see \cite[17.8(1)]{Brzezinski/Wisbauer:2003}).
\begin{multline*}
\begin{split}
H(b.f) & =H(i_R(b)*^{r}f) &
=\bigg[\Big(\phi_{\Lambda,B}(f\delta_{\Lambda})
\phi_{\Lambda,B}\big(i_R(b)\delta_{\Lambda}\big)
[\xymatrix{\Lambda \ar[r]^-\simeq &
{B\tensor{B}{\Lambda}}}]\Big)\cotensor{\coring{C}}
X\bigg]\xi_B(1_B).
\end{split}
\end{multline*}
On the other hand,
\begin{multline*}
\begin{split}
\phi_{\Lambda,B}(i_R(b)\delta_{\Lambda})\cotensor
{\coring{C}}X & =
[\xymatrix{B\tensor{B}{\Lambda}\cotensor{\coring{C}}X\ar[r]^-\simeq
&
{\Lambda\cotensor{\coring{C}}X}}]\big(\epsilon(b-)\tensor{B}{\Lambda}
\cotensor{\coring{C}}X\big)
\Delta \\
& =b1_{\Lambda\cotensor{\coring{C}}X}.
\end{split}
\end{multline*}
Hence $H(b.f)=bH(f)$.
\end{proof}

\medskip

The last result is a generalization of \cite[Theorem 3.2(1), Theorem
3.5(1)]{Brzezinski/Gomez:2003}.

\begin{corollary}\label{comatrix}
Let $_BM_A$ be a $(B,A)$-bimodule such that $M_A$ is finitely
generated projective. Then \begin{enumerate}[(1)]
\item $_AM^*_B$ is
a separable bimodule if and only if the comatrix coring
$M^*\tensor{B}M$ is a cosplit $A$-coring.
\item If $M$ is a separable bimodule, then the
comatrix coring $M^*\tensor{B}M$ is a coseparable $A$-coring. \item
If $M$ is a Frobenius bimodule, then the comatrix coring
$M^*\tensor{B}M$ is a Frobenius $A$-coring.
\end{enumerate}
\end{corollary}

\begin{proof}
 Define $\psi:B\rightarrow M\tensor{A}M^*$ by
$b\mapsto
 \sum_{i=1}^{n}be_i\otimes
e_i^*=\sum_{i=1}^{n}e_i\otimes e_i^*b$, where
$\{e_i,e_i^*\}_{i\in\{1,\dots,n\}}$ is a dual basis basis of $M_A$,
and $\omega:M^*\tensor{B}M\rightarrow A$ by $\varphi\otimes
m\mapsto\varphi(m)$ (the evaluation map). Then
$(A,B,{_AM^*_B},{_BM_A},\omega,\psi)$ is a comatrix coring context
(see \cite[p. 3]{Brzezinski/Gomez:2003}). From \cite[Proposition
2.8]{Zarouali:2004}, we have the adjunction
$(-\tensor{B}{M},-\tensor{A}{M^*})$, in particular,
$\Lambda=M^*\in{_A\mathcal{M}_B}$ is quasi-finite as a right
$B$-(co)module, and the cohom functor is $-\tensor{B}{M}$.

(1) It is clear from the definition of a separable bimodule (see
\cite{Sugano:1971} or \cite{Kadison:1999}), that $_AM^*_B$ is a
separable bimodule if and only if $\omega$ is an $A$-bimodule
retraction.

(2) Let $S=\operatorname{End}_A(M_A)$. Define $\lambda:B\rightarrow
S$ by $b\mapsto\lambda_B:M\rightarrow M,[m\mapsto bm]$, and
$\phi:S\rightarrow M\tensor{A}{M^*}$ by $s\mapsto
\sum_is(e_i)\otimes e_i^*$, where
$\{e_i,e_i^*\}_{i\in\{1,\dots,n\}}$ is a dual basis of $M_A$. $\phi$
is an isomorphism of $B$-bimodule, with inverse map $m\otimes
\varphi\mapsto[x\mapsto m\varphi(x)]$. It is easy to see that
$\psi=\phi\circ\lambda$. From Sugano's theorem \cite[Theorem 1,
Proposition 2]{Sugano:1971} (see also \cite[Theorem
3.1]{Kadison:1999}), the bimodule $M$ is separable if and only if
$\lambda$ is a split extention, i.e., $\lambda$ is a $B$-bimodule
section. However that is equivalent to $\psi$ is a $B$-bimodule
section.

(3) This is \cite[Theorem 3.7(1)]{Brzezinski/Gomez:2003}.
\end{proof}

\section{Applications to equivalences of categories of
comodules}\label{equivalences}

In this section we will generalize and improve the main results
concerning equivalences between categories of comodules given in
\cite{Takeuchi:1977}, \cite{AlTakhman:2002} and
\cite{Brzezinski/Wisbauer:2003}. We also give new characterizations
of equivalences between categories of comodules over separable
corings or corings with a duality.

\medskip

A functor $S:\cat{C}\rightarrow \cat{D}$ is called an equivalence,
if there exists a functor $T:\cat{D}\rightarrow \cat{C}$ with
natural isomorphisms $\eta:1\rightarrow TS$ and
$\varepsilon:ST\rightarrow 1$. We recall from \cite[p.
93]{MacLane:1998} that an adjoint equivalence of categories is an
adjunction $(S,T,\eta,\varepsilon)$ in which both the unit
$\eta:1\rightarrow TS$ and the counit $\varepsilon:ST\rightarrow 1$
are natural isomorphisms. We have, $(S,T)$ is a pair of inverse
equivalences if and only if $S$ and $T$ are part of an adjoint
equivalence $(S,T,\eta,\varepsilon)$. The proof is analogous to that
of \cite[Theorem IV.4.1]{MacLane:1998}: The ``if'' part is trivial.
For the ``only if'' part, let $\eta:1\rightarrow TS$ be a natural
isomorphism, then $\varphi_{C,D}(\alpha)= T(\alpha).\eta_{C}$ for
$\alpha:S(C)\rightarrow D$, is a natural transformation in $C$ and
$D$. Since $T$ is faithful and full, $\varphi$ is a natural
isomorphism, and $T$ is a right adjoint to $S$. Therefore the counit
of this adjunction $\varepsilon:ST\rightarrow 1$ satisfies
$T(\varepsilon_{D}).\eta_{T(D)}=1_{T(D)}$ for every $D\in\cat{D}$,
and consequently $T(\varepsilon_{D})=(\eta_{T(D)})^{-1}$ is
invertible. Since $T$ is faithful and full, $\varepsilon_{D}$ is
also invertible, and $\varepsilon$ is a natural isomorphism (to show
that $\varepsilon$ is a natural isomorphism we can also use
\cite[Theorem IV.3.1]{MacLane:1998}).

\medskip

The following result generalizes the case of the category of modules
\cite[II (2.4)]{Bass:1968}.

\begin{proposition}\label{equivalence1}
Suppose that $_A\coring{C}$, $\coring{C}_A$, $_B\coring{D}$ and
$\coring{D}_B$ are flat. Let
$X\in{}^{\coring{C}}\mathcal{M}^\coring{D}$ and
$\Lambda\in{}^{\coring{D}}\mathcal{M}^\coring{C}$.
\item The
following statements are equivalent:
\begin{enumerate}[(1)]
\item $(-\square
_{\coring{C}}X,-\square_{\coring{D}}\Lambda)$ is a pair of inverse
equivalences;
\item there exist bicomodule isomorphisms
\[
 f:X\square_{\coring{D}}\Lambda
\rightarrow\coring{C}\quad \text{and}\quad
g:\Lambda\square_{\coring{C}}X\rightarrow\coring{D}
\]
in $^\coring{C}\mathcal{M}\coring{C}$ and
$^{\coring{D}}%
\mathcal{M}^\coring{D}$ respectively, such that
\begin{enumerate}[(a)] \item $_A X$ and $_B \Lambda$
are flat, and
$\omega_{X,\Lambda}=\rho{}_X\otimes_B\Lambda-X\otimes_A\rho{}_\Lambda$
is pure in $\lmod{A}$ and
$\omega_{\Lambda,X}=\rho{}_\Lambda\otimes_AX-\Lambda\otimes_B\rho{}_X$
is pure in $\lmod{B}$, or \item ${}_{\coring{C}}X$ and
${}_{\coring{D}}\Lambda$ are coflat.
\end{enumerate}
\end{enumerate}
In such a case the diagrams
\begin{equation}
\xymatrix{\Lambda\cotensor{\coring{C}}X\cotensor{\coring{D}}\Lambda
\ar[rr]^{\Lambda\cotensor{\coring{C}}f}
\ar[d]_{g\cotensor{\coring{D}}\Lambda} &&
\Lambda\cotensor{\coring{C}}\coring{C} \ar[d]^\simeq
\\\coring{D}\cotensor{\coring{D}}\Lambda
\ar[rr]^\simeq && \Lambda}
\xymatrix{X\cotensor{\coring{D}}\Lambda\cotensor{\coring{C}}X
\ar[rr]^{f\cotensor{\coring{C}}X}
\ar[d]_{X\cotensor{\coring{D}}g}&&
\coring{C}\cotensor{\coring{C}}X\ar[d]^\simeq
\\X\cotensor{\coring{D}}\coring{D} \ar[rr]^\simeq &&
X}
\end{equation}
commute.

If $A$ and $B$ are von Neumann regular rings, or if $\coring{C}$ and
$\coring{D}$ are coseparable corings (without $\coring{C}_A$ and
$\coring{D}_B$ are flat), the conditions (a) and (b) can be deleted.
\end{proposition}

\begin{proof}
Clear from \cite[Lemma 2.6, Proposition 2.7]{Zarouali:2004} and the
above mentioned consideration.
\end{proof}

\medskip

Following \cite{AlTakhman:2002} and \cite{Brzezinski/Wisbauer:2003},
we state the following definition.

\begin{definition}\label{7}
A bicomodule $N\in{}^{\coring{C}}\mathcal{M}^\coring{D}$ is called
an \emph{injector} (resp. \emph{a cogenerator preserving}, resp.
\emph{an injector-cogenerator}) as a right $\coring{D}$-comodule if
the functor
$-\otimes_AN:\mathcal{M}^{A}\rightarrow\mathcal{M}^\coring{D}$
preserves injective (resp. cogenerator, resp. injective cogenerator)
objects.
\end{definition}

\medskip

The version for module categories of the following lemma is given in
\cite[Exercise 20.8]{Anderson/Fuller:1992}.

\begin{lemma}\label{injector-cogenerator}
Let $\cat{C}$ and $\cat{D}$ be two abelian categories, and let
$S:\cat{C}\rightarrow \cat{D}$ be a functor. Let
$T:\cat{D}\rightarrow \cat{C}$ be a right adjoint of $S$. We have
the following properties:
\begin{enumerate}[(1)] \item If $N\in\cat{D}$ is
injective and $S$ is exact, then $T(N)\in\cat{C}$ is injective.
\item If $N\in\cat{D}$ is cogenerator and $T(N)\in\cat{C}$ is
injective, then $S$ is exact.
\item If $N\in\cat{D}$ is cogenerator and $S$ is
faithful, then $T(N)\in\cat{C}$ is a cogenerator. \item If
$T(N)\in\cat{C}$ is a cogenerator for some $N\in\cat{D}$, then $S$
is faithful.

If moreover $\cat{D}$ is a category with sufficiently many
injectives (e.g. if $\cat{D}$ is a Grothendieck category), then the
following property holds: \item $S$ is exact if and only if $T$
preserves injective objects.

If moreover $\cat{D}$ is an AB 3 category with generators and
sufficiently many injectives (e.g. if $\cat{D}$ is a Grothendieck
category), then the following properties hold: \item $S$ is faithful
if and only if $T$ preserves cogenerator objects.\item $S$ is
faithfully exact if and only if $T$ preserves injective cogenerator
objects.
\end{enumerate}

In particular, let $N\in{}^\coring{C}\mathcal{M}^\coring{D}$ be a
bicomodule, quasi-finite as a right $\coring{D}$-comodule, such that
$_A\coring{C}$ and $_B\coring{D}$ are flat. $N$ is an injector
(resp. a cogenerator preserving, resp. an injector-cogenerator) as a
right $\coring{D}$-comodule if and only if the cohom functor
$h_{\coring{D}}(N,-)$ is exact (resp. faithful, resp. faithfully
exact).
\end{lemma}

\begin{proof}
(1) is \cite[Lemma 3.2.7]{Popescu:1973}. (2), (3) and (4) are clair
from the theorem of Freyd \cite[Theorem II.7.1]{Mitchell:1965}. (5)
is \cite[Theorem 3.2.8]{Popescu:1973}. $(6)$ is an immediate
consequence of (3) and (4) and the fact that every Grothendieck
category has a cogenerator object (see \cite[Lemma
3.7.12]{Popescu:1973}). The ``only if'' part of $(7)$ is obvious
from (1) and (3). The ``if'' part is an immediate consequence of (2)
and (4) and the fact that every Grothendieck category has an
injective cogenerator object (see \cite[Lemma
3.7.12]{Popescu:1973}).
\end{proof}

\medskip

The first part of the following is
\cite[23.10]{Brzezinski/Wisbauer:2003}. We think that the proof we
give here is more clear.

\begin{proposition}\label{equivalence2}
Let $_A\coring{C}$ be flat and let
$\Lambda\in{}_B\mathcal{M}^\coring{C}$ be quasi-finite as a right
$\coring{C}$-comodule. Let $\cohom{\coring{C}}{\Lambda}{-}$ be the
cohom functor of $\Lambda$. Denote by $\e{\coring{C}}{\Lambda}$ the
coendomorphism coring of $\Lambda$. Suppose that
$\e{\coring{C}}{\Lambda}$ is flat as left $B$-module.
\begin{enumerate}[(1)] \item If $\coring{C}_A$ is flat
and $\e{\coring{C}}{\Lambda}$ is flat as right $B$-module, and
$\Lambda$ is an injector-cogenerator as a right
$\coring{C}$-comodule, then the functors
\begin{equation}\label{functors}
-\square_{\e{\coring{C}}{\Lambda}}\Lambda:\mathcal{M}^{\e{\coring{C}}{\Lambda}}
\rightarrow \mathcal{M}^\coring{C},\qquad
\cohom{\coring{C}}{\Lambda}{-}:\mathcal{M}^\coring{C}\rightarrow
\mathcal{M}^{\e{\coring{C}}{\Lambda}},
\end{equation}
are inverse equivalences. \item If $\coring{C}$ and
$\e{\coring{C}}{\Lambda}$ are coseparable corings, and $\Lambda$ is
a cogenerator preserving as a right $\coring{C}$-comodule, then the
functors of \eqref{functors} are inverse equivalences.
\end{enumerate}
\end{proposition}

\begin{proof}
We will prove at the same time the two statements. Let
$\theta:1_{\mathcal{M}^\coring{C}}\rightarrow
\cohom{\coring{C}}{\Lambda}{-}\tensor{B}{\Lambda}$ be the unit of
the adjunction
$(\cohom{\coring{C}}{\Lambda}{-},-\tensor{B}{\Lambda})$. By
\cite[Theorem 2.3]{Zarouali:2004},
$\cohom{\coring{C}}{\Lambda}{-}\simeq
-\square_{\coring{C}}\cohom{\coring{C}}{\Lambda}{\coring{C}}$. In
the case (1), $\cohom{\coring{C}}{\Lambda}{\coring{C}}$ is coflat as
a left $\coring{C}$-comodule. We put
$\coring{D}=\e{\coring{C}}{\Lambda}$. By \cite[Proposition
3.4]{Gomez:2002},
$\delta_{\Lambda}:\coring{D}\rightarrow\Lambda\square_{\coring{C}}\cohom{\coring{C}}
{\Lambda}{\coring{C}}$ is $\coring{D}-\coring{D}$-bicolinear.
Therefore, $-\square_{\coring{D}}\Lambda\square_{\coring{C}}
\cohom{\coring{C}}{\Lambda}{\coring{C}}\simeq-\square_{\coring{D}}\coring{D}$
is exact. Hence, in the case (1), the functor
$-\square_{\coring{D}}\Lambda$ is exact (since
$\cohom{\coring{C}}{\Lambda}{-}$ is faithful). From the proof of
Proposition \ref{comatrixproperty1}(1), the unit of the adjunction
$(-\square_{\coring{C}}\cohom{\coring{C}}{\Lambda}{\coring{C}},
-\square_{\coring{D}}\Lambda)$
is
$$\xymatrix@1{\eta
:1_{\mathcal{M}^\coring{C}}\ar[r]^\simeq & -\square
_{\coring{C}}\coring{C}
\ar[rr]^-{-\square_{\coring{C}}\theta_{\coring{C}}}
& &
-\square_{\coring{C}}\cohom{\coring{C}}{\Lambda}{\coring{C}}
\square_{\coring{D}}\Lambda}.$$
From Proposition \ref{delta}, $\delta_{\Lambda}$ is the unique map
making the following diagram commutative
\[
\xymatrix{\Lambda \ar[d]_\simeq\ar[rr]^\simeq & &
\Lambda\square_{\coring{C}}\coring{C}\ar[d]^{\Lambda
\square_{\coring{C}}\theta_{\coring{C}}}
\\
\coring{D}\square_{\coring{D}}\Lambda
\ar[rr]_-{\delta_{\Lambda}\square_{\coring{D}}\Lambda} & &
\Lambda\square_{\coring{C}}(\cohom{\coring{C}}{\Lambda}
{\coring{C}}\square_{\coring{D}}\Lambda).}
\]
By \cite[Lemma 2.6, Poposition 2.7]{Zarouali:2004}, the counit of
the adjunction
$(-\square_{\coring{C}}\cohom{\coring{C}}{\Lambda}{\coring{C}},
-\square_{\coring{D}}\Lambda)$
is
$$\xymatrix@1{\varepsilon:-\square_{\coring{D}}
\Lambda\square_{\coring{C}}
\cohom{\coring{C}}{\Lambda}{\coring{C}}
\ar[rr]^-{-\square_{\coring{D}}\delta^{-1}_{\Lambda}} & &
-\square_{\coring{D}}\coring{D}\ar[r]^\simeq &
1_{\mathcal{M}^\coring{D}},} $$ and we have the commutative diagram
\[
\xymatrix{\cohom{\coring{C}}{\Lambda}{\coring{C}}
\ar[d]_\simeq\ar[rr]^\simeq
& &
\coring{C}\square_{\coring{C}}\cohom{\coring{C}}{\Lambda}{\coring{C}}
\ar[d]^{\theta_{\coring{C}}\square_{\coring{C}}\cohom{\coring{C}}
{\Lambda}{\coring{C}}}
\\
\cohom{\coring{C}}{\Lambda}{\coring{C}}\square_{\coring{D}}\coring{D}
& &
\cohom{\coring{C}}{\Lambda}{\coring{C}}\square_{\coring{D}}
\Lambda\square_{\coring{C}}
\cohom{\coring{C}}{\Lambda}{\coring{C}}
\ar[ll]_-{\cohom{\coring{C}}{\Lambda}{\coring{C}}\square_{\coring{D}}
\delta^{-1}_{\Lambda}}.}
\]
Then,
$\theta_{\coring{C}}\square_{\coring{C}}
\cohom{\coring{C}}{\Lambda}{\coring{C}}$
is an isomorphism. Since
$-\square_{\coring{C}}\cohom{\coring{C}}{\Lambda}{\coring{C}}$ is
faithful, $\theta_{\coring{C}}$ is also an isomorphism. Finally, the
unit and the counit of the adjunction
$(-\square_{\coring{C}}\cohom{\coring{C}}{\Lambda}{\coring{C}},
-\square_{\coring{D}}\Lambda)$
are natural isomorphisms.
\end{proof}

\medskip

The first part of the following is contained in \cite[Theorem
3.10]{ElKaoutit/Gomez:2003}.

\begin{corollary}
Let ${}_BM_A$ be a $(B,A)$-bimodule such that $M_A$ is finitely
generated projective. \begin{enumerate}[(1)] \item The following
statements are equivalent \begin{enumerate}[(a)]\item
${}_A(M^*\otimes_BM)$ is flat and
$-\otimes_BM:\mathcal{M}_B\rightarrow\mathcal{M}^{M^*\otimes_BM}$ is
an equivalence of categories; \item $_B{}M$ is faithfully flat.
\end{enumerate}
\item If $\mathcal{M}^{M^*\otimes_BM}$ is an abelian
category, and $M^*\otimes_BM$ is an $A$-coseparable coring, then the
following statements are equivalent\begin{enumerate}[(a)]\item
$-\otimes_BM:\mathcal{M}_B\rightarrow\mathcal{M}^{M^*\otimes_BM}$ is
an equivalence of categories; \item $_BM$ is completely faithful.
\end{enumerate}
\end{enumerate}
\end{corollary}

\begin{proof}
It suffices to take $\Lambda=M^*\in{}_A\mathcal{M}_B$ in Proposition
\ref{equivalence2}.
\end{proof}

The first part of the following is contained in
\cite[23.12]{Brzezinski/Wisbauer:2003}. We think that the proof we
give here is more clear.

\begin{proposition}\label{equivalence3}
Suppose that $_A\coring{C}$ and $_B\coring{D}$ are flat, and let
$X\in{}^{\coring{C}}\mathcal{M}^\coring{D}$ and
$\Lambda\in{}^{\coring{D}}\mathcal{M}^\coring{C}$.

If $\coring{C}_A$ and $\coring{D}_B$ are flat, then the following
statements are equivalent
\begin{enumerate}[(1)]
\item $(-\square
_{\coring{C}}X,-\square_{\coring{D}}\Lambda)$ is a pair of inverse
equivalences; \item $\Lambda$ is quasi-finite injector-cogenerator
as a right
$\coring{C}%
$-comodule, $\e{\coring{C}}{\Lambda}\simeq \coring{D}$ as corings
and $X\simeq \cohom{\coring{C}}{\Lambda}{\coring{C}}$ in
$^{\coring{C}}\mathcal{M}^\coring{D}$;\item $X$ is quasi-finite
injector-cogenerator as a right $\coring{D}%
$-comodule, $\e{\coring{D}}{X}\simeq \coring{C}$ as corings and
$\Lambda\simeq \cohom{\coring{D}}{X}{\coring{D}}$ in
$^{\coring{D}}\mathcal{M}^\coring{C}$.

If moreover $\coring{C}$ and $\coring{D}$ are coseparable or
cosemisimple, then (1) is equivalent to
\item $\Lambda$ is quasi-finite cogenerator preserving
as a right
$\coring{C}%
$-comodule, $\e{\coring{C}}{\Lambda}\simeq \coring{D}$ as corings
and $X\simeq \cohom{\coring{C}}{\Lambda}{\coring{C}}$ in
$^{\coring{C}}\mathcal{M}^\coring{D}$;\item $X$ is quasi-finite
cogenerator preserving as a right $\coring{D}%
$-comodule, $\e{\coring{D}}{X}\simeq \coring{C}$ as corings and
$\Lambda\simeq \cohom{\coring{D}}{X}{\coring{D}}$ in
$^{\coring{D}}\mathcal{M}^\coring{C}$.
\end{enumerate}
\end{proposition}

\begin{proof}
At first we will prove the first part.

$(1)\Rightarrow(2)$ By \cite[Proposition 2.9]{Zarouali:2004},
$\Lambda$ is quasi-finite injector as a right
$\coring{C}%
$-comodule, and $X\simeq \cohom{\coring{C}}{\Lambda}{\coring{C}}$ in
$^{\coring{C}}\mathcal{M}^\coring{D}$. Therefore, $-\square
_{\coring{C}}X\simeq \cohom{\coring{C}}{\Lambda}{-}$ is faithful. By
Lemma \ref{injector-cogenerator}, $\Lambda$ is quasi-finite
injector-cogenerator as a right $\coring{C}$-comodule. Finally, by
Proposition \ref{comatrixproperty1}, $\e{\coring{C}}{\Lambda}\simeq
\coring{D}$ as corings.

$(2)\Rightarrow(1)$ By \cite[Proposition 2.9]{Zarouali:2004},
$-\square _{\coring{C}}X\simeq \cohom{\coring{C}}{\Lambda}{-}$.
Hence (1) follows obviously from Proposition \ref{equivalence2}(1).

$(1)\Leftrightarrow(3)$ Follows by symmetry.

 Now we will prove the second part. The case ``cosemisimple''
is obvious from the first part. It suffices to show the
``coseparable'' case.

$(1)\Rightarrow(4)$ Obvious from the first part.

$(4)\Rightarrow(1)$ By \cite[Proposition 2.7]{Zarouali:2004},
$-\square _{\coring{C}}X\simeq \cohom{\coring{C}}{\Lambda}{-}$.
Hence (1) follows obviously from Proposition \ref{equivalence2}(2).

$(1)\Leftrightarrow(5)$ Follows by symmetry.
\end{proof}

\medskip

As an immediate consequence of \cite[Proposition 2.7]{Zarouali:2004}
and Proposition \ref{equivalence3}, we get the two following
theorems.

\begin{theorem}\label{equivalence4}
Suppose that $_A\coring{C}$, $\coring{C}_A$, $_B\coring{D}$ and
$\coring{D}_B$ are flat, and let
$X\in{}^{\coring{C}}\mathcal{M}%
^{\coring{D}}$ and
$\Lambda\in{}^{\coring{D}}\mathcal{M}^\coring{C}$. The following
statements are equivalent:
\begin{enumerate}[(1)] \item $(-\square
_{\coring{C}}X,-\square_{\coring{D}}\Lambda)$ is a pair of inverse
equivalences with $X_{\coring{D}}$ and $\Lambda_{\coring{C}}$ are
coflat; \item $(\Lambda \cotensor{\coring{C}}-,X
\cotensor{\coring{D}}-)$ is a pair of inverse equivalences with
$_{\coring{C}}X$ and $_{\coring{D}}\Lambda$ are coflat; \item
$\Lambda$ is quasi-finite injector-cogenerator as a right
$\coring{C}%
$-comodule with $X_{\coring{D}}$ and $\Lambda_{\coring{C}}$ are
coflat, $\e{\coring{C}}{\Lambda}\simeq \coring{D}$ as corings and
$X\simeq \cohom{\coring{C}}{\Lambda}{\coring{C}}$ in
$^{\coring{C}}\mathcal{M}^\coring{D}$;\item $X$ is quasi-finite
injector-cogenerator as a right $\coring{D}%
$-comodule with $X_{\coring{D}}$ and $\Lambda_{\coring{C}}$ are
coflat, $\e{\coring{D}}{X}\simeq \coring{C}$ as corings and
$\Lambda\simeq \cohom{\coring{D}}{X}{\coring{D}}$ in
$^{\coring{D}}\mathcal{M}^\coring{C}$;\item $\Lambda$ is
quasi-finite injector-cogenerator coflat on both sides,
$\e{\coring{C}}{\Lambda}\simeq \coring{D}$ as corings, and $X\simeq
\cohom{\coring{C}}{\Lambda}{\coring{C}}$ in
$^{\coring{C}}\mathcal{M}^\coring{D}$;\item $X$ is quasi-finite
injector-cogenerator coflat on both sides, $\e{\coring{D}}{X}\simeq
\coring{C}$ as corings, and $\Lambda\simeq
\cohom{\coring{D}}{X}{\coring{D}}$ in
$^{\coring{D}}\mathcal{M}^\coring{C}.$
\end{enumerate}
\end{theorem}

\medskip

\begin{theorem}\label{equivalence5}
Suppose that $_A\coring{C}$, $\coring{C}_A$, $_B\coring{D}$ and
$\coring{D}_B$ are flat, and let
$X\in{}^{\coring{C}}\mathcal{M}%
^{\coring{D}}$ and
$\Lambda\in{}^{\coring{D}}\mathcal{M}^\coring{C}$. If $\coring{C}$
and $\coring{D}$ are coseparable or cosemisimple (resp. $A$ and $B$
are von Neumann regular ring), then the following statements are
equivalent:
\begin{enumerate}[(1)] \item $(-\square
_{\coring{C}}X,-\square_{\coring{D}}\Lambda)$ is a pair of inverse
equivalences; \item $(\Lambda \cotensor{\coring{C}}-,X
\cotensor{\coring{D}}-)$ is a pair of inverse equivalences; \item
$\Lambda$ is quasi-finite cogenerator preserving (resp.
injector-cogenerator) as a right
$\coring{C}%
$-comodule, $\e{\coring{C}}{\Lambda}\simeq \coring{D}$ as corings
and $X\simeq \cohom{\coring{C}}{\Lambda}{\coring{C}}$ in
$^{\coring{C}}\mathcal{M}^\coring{D}$;\item $X$ is quasi-finite
cogenerator preserving (resp. injector-cogenerator) as a right
$\coring{D}%
$-comodule, $\e{\coring{D}}{X}\simeq \coring{C}$ as corings and
$\Lambda\simeq \cohom{\coring{D}}{X}{\coring{D}}$ in
$^{\coring{D}}\mathcal{M}^\coring{C}$;\item $\Lambda$ is
quasi-finite cogenerator preserving (resp. injector-cogenerator) on
both sides, $\e{\coring{C}}{\Lambda}\simeq \coring{D}$ as corings,
and $X\simeq \cohom{\coring{C}}{\Lambda}{\coring{C}}$ in
$^{\coring{C}}\mathcal{M}^\coring{D}$;\item $X$ is quasi-finite
cogenerator preserving (resp. injector-cogenerator) on both sides,
$\e{\coring{D}}{X}\simeq \coring{C}$ as corings, and $\Lambda\simeq
\cohom{\coring{D}}{X}{\coring{D}}$ in
$^{\coring{D}}\mathcal{M}^\coring{C}.$
\end{enumerate}
\end{theorem}

\begin{remark}\label{Morita}
As a relevant consequence of the last theorem, we can prove the
Morita's characterization of equivalence \cite[Theorem
22.2]{Anderson/Fuller:1992}: Let $A$  and $B$ be rings and let
$$F:\rmod{A}\rightarrow \rmod{B}\qquad\textrm{and
}\qquad G:\rmod{B}\rightarrow \rmod{A}$$ be additive functors. Then
the following statements are equivalent
\begin{enumerate}[(1)]\item $F$ and $G$ are inverse
equivalences;
\item there exists a bimodule $_AM_B$ such that:
\begin{enumerate}[(a)]\item $_AM$ and $M_B$ are
progenerators (i.e., finitely generated, projective and generators),
\item $_AM_B$
is (faithfully) balanced (see \cite[p. 60]{Anderson/Fuller:1992}),
\item $F\simeq -\tensor{A}{M}\qquad\textrm{and}\qquad
G\simeq \hom{B}{M}{-}$;
\end{enumerate}
\item there exists a bimodule $_AM_B$ such that: $M_B$
is finitely generated projective, $_AM$ is completely faithful, the
evaluation map $M^*\tensor{A}{M}\rightarrow B$ is an isomorphism,
and
\\ $F\simeq -\tensor{A}{M}\qquad\textrm{and}\qquad
G\simeq \hom{B}{M}{-}$.
\end{enumerate}
Moreover in such a case, $_BM^*$ and $M^*_A$ are progenerators, and
$$F\simeq \hom{A}{M^*}{-}\qquad\textrm{and}\qquad
G\simeq -\tensor{B}{M^*}.$$
\end{remark}

\medskip
Now, we will study when the category $\rcomod{\coring{C}}$ for some
coring $\coring{C}$ is equivalent to a category of modules
$\rmod{B}$.

Let $\Sigma\in {}_B\rcomod{\coring{C}}$ such that $\Sigma_A$ is
finitely generated projective with a dual basis $\{e_i,e_i^*\}_i$,
and let $M\in \rcomod{\coring{C}}$. By \cite[Proposition
1.4]{Caenepeel/DeGroot/Vercruysse:unp}, the canonical isomorphism
$\hom{A}{\Sigma}{M}\to M\tensor{A}\rdual{\Sigma}$ yields an
isomorphism $\hom{\coring{C}}{\Sigma}{M}\to
M\cotensor{\coring{C}}\rdual{\Sigma}$, where the left coaction on
$\Sigma^*$ is given by
$$ \lambda_{\Sigma^*}:\Sigma^*\to
\coring{C}\tensor{A}\Sigma^*,\; \lambda_{\Sigma^*}(f)=\sum_i
f(e_{i(0)})e_{i(1)}\tensor{A}e_i^*.$$

Moreover, we have an adjoint pair $(F,G)$, where
$$F=-\tensor{B}\Sigma:\rmod{B}\to
\rcomod{\coring{C}},$$ and
$$G=\hom{\coring{C}}{\Sigma}{-}:\rcomod{\coring{C}}\to
\rmod{B}.$$ The unit and the counit of this adjunction are given by:
For $N\in\rmod{B}$,
$$\upsilon_N:N\to
\hom{\coring{C}}{\Sigma}{N\tensor{B}\Sigma},\;
\upsilon_N(n)(u)=n\tensor{B}u,$$ or
$$\upsilon_N:N\to
(N\tensor{B}\Sigma)\cotensor{\coring{C}}\rdual{\Sigma},\;
\upsilon_N(n)=\sum_i(n\tensor{B}e_i)\tensor{A}e_i^*,$$ and for $M\in
\rcomod{\coring{C}}$:
$$\zeta_M:\hom{\coring{C}}{\Sigma}{M}\tensor{B}\Sigma\to
M,\;\zeta_M(\varphi\tensor{B}u)=\varphi(u),$$ or
$$\zeta_M:(M\cotensor{\coring{C}}\rdual{\Sigma})\tensor{B}\Sigma\to
M, \;
\zeta_M\Big((\sum_jm_j\tensor{A}f_j)\tensor{B}u\Big)=\sum_jm_jf_j(u)$$
(see \cite[Proposition 1.5]{Caenepeel/DeGroot/Vercruysse:unp}).

Let us define the map
$$\mathbf{can}:\rdual{\Sigma}\tensor{B}\Sigma\to
\coring{C},\; \mathbf{can}(f\tensor{B}u)=f(u_{(0)})u_{(1)}.$$ From
\cite[Lemma 3.1]{Caenepeel/DeGroot/Vercruysse:unp}, we have
$\mathbf{can}$ is a morphism of corings. It follows from this that
$\mathbf{can}$ is a $\coring{C}$-bicolinear map. Moreover, we can
verify easily that for every $M\in \rcomod{\coring{C}}$, the
following diagram is commutative:
\begin{equation}\label{counitmod}
\xymatrix{(M\cotensor{\coring{C}}\rdual{\Sigma})\tensor{B}\Sigma
\ar[rr]^-{\zeta_M} \ar[d]^{\psi_M} & & M\ar[d]^\simeq
\\
M\cotensor{\coring{C}}(\rdual{\Sigma}\tensor{B}\Sigma)\ar[rr]^-
{M\cotensor{\coring{C}}\mathbf{can}}
& & M\cotensor{\coring{C}}\coring{C},}
\end{equation}
where $\psi_M$ is the canonical map. We obtain that if $\psi_M$ is
isomorphism, for every $M\in \rcomod{\coring{C}}$ (for example if
$_B\Sigma$ is flat or if $\coring{C}$ is coseparable), then $G$ is
fully faithful if and only if $\mathbf{can}$ is an isomorphism (see
\cite[Theorem IV.3.1]{MacLane:1998}).

Let us consider the map $\upsilon:B\to
\Sigma\cotensor{\coring{C}}\rdual{\Sigma},\;
\upsilon(b)=\sum_ibe_i\tensor{A}e_i^*$ $(b\in B).$ For every $N\in
\rmod{B}$, we have the following commutative diagram
\begin{equation}\label{unitmod}
\xymatrix{N \ar[rr]^-{\upsilon_N}\ar[d]_\simeq & &
(N\tensor{B}\Sigma)\cotensor{\coring{C}}\rdual{\Sigma}
\\ N\tensor{B}B\ar[rr]^-{N\tensor{B}\upsilon}
& &
N\tensor{B}(\Sigma\cotensor{\coring{C}}\rdual{\Sigma})\ar[u]_{\psi_N},}
\end{equation}
where $\psi_N$ is the canonical map. Let $\phi:B\to
\operatorname{End}_\coring{C}(\Sigma)$ be the canonical morphism of
$k$-algebras which define the left action on $\Sigma$
($\phi(b)(u):=bu$) (see Section 1 of \cite{Zarouali:2004}). We have
$\upsilon:\xymatrix{B\ar[r]^-\phi &
\operatorname{End}_\coring{C}(\Sigma)\ar[r]^-\simeq &
\Sigma\cotensor{\coring{C}}\rdual{\Sigma}}$. We obtain that if
$\psi_N$ is an isomorphism, for every $N\in \rmod{B}$ (for example
if $\Sigma_\coring{C}$ is projective or if $\coring{C}$ is
coseparable), then $F$ is fully faithful if and only if $\phi$ is an
isomorphism (see \cite[Theorem IV.3.1]{MacLane:1998}).

Furthermore, if $_B\Sigma$ is flat or $\coring{C}$ is coseparable,
then we have the commutative diagram:
\begin{equation}\label{unitcounitmod}
\xymatrix{\Sigma\ar[rr]^\simeq \ar[d]_\simeq & &
B\tensor{B}\Sigma\ar[d]^{\upsilon\tensor{B}\Sigma}\\
\Sigma\cotensor{\coring{C}}\coring{C} & & \Sigma
\cotensor{\coring{C}}\Sigma^*\tensor{B}\Sigma
\ar[ll]^{\Sigma\cotensor{\coring{C}}\mathbf{can}}.}
\end{equation}

Now we are ready to state and prove the following theorem. In the
particular case where $B=\operatorname{End}_\coring{C}(\Sigma)$, the
first part of Theorem is known, see \cite[Theorem
3.2]{ElKaoutit/Gomez:2003} and
\cite[18.27]{Brzezinski/Wisbauer:2003}.

\begin{theorem}\label{equivmodcomod}
Let $\Sigma\in {}_B\rcomod{\coring{C}}$ such that $\Sigma_A$ is
finitely generated projective. Let $F=-\tensor{B}\Sigma:\rmod{B}\to
\rcomod{\coring{C}}$ and
$G=\hom{\coring{C}}{\Sigma}{-}:\rcomod{\coring{C}}\to \rmod{B}$. The
following statements are equivalent:
\begin{enumerate}[(1)]
\item$(F,G)$ is a pair of inverse equivalences with
$_A\coring{C}$ flat;
\item $_B\Sigma$ is flat, $\Sigma_\coring{C}$ is
projective, and $\mathbf{can}$, and the morphism of $k$-algebras
$\phi:B\to \operatorname{End}_\coring{C}(\Sigma)$ defined by
$\phi(b)(u)=bu$ for $b\in B,u\in \Sigma$, are isomorphisms;
\item $_A\coring{C}$ is flat, $\Sigma_\coring{C}$ is
a projective generator, and the morphism of $k$-algebras $\phi:B\to
\operatorname{End}_\coring{C}(\Sigma)$ defined by $\phi(b)(u)=bu$
for $b\in B,u\in \Sigma$, is an isomorphism;
\item $_B\Sigma$ is faithfully flat, and
$\mathbf{can}$ is an isomorphism.
\end{enumerate}
If moreover $\coring{C}$ is coseparable, and $_A\coring{C}$ is
projective, then the following statements are equivalent:
\begin{enumerate}[(1)]\item $(F,G)$ is a pair of
inverse equivalences;
\item $\Sigma_\coring{C}$ is a generator, and the
morphism of $k$-algebras $\phi:B\to
\operatorname{End}_\coring{C}(\Sigma)$ defined by $\phi(b)(u)=bu$
for $b\in B,u\in \Sigma$, is an isomorphism;
\item the morphism of $k$-algebras $\phi:B\to
\operatorname{End}_\coring{C}(\Sigma)$ defined by $\phi(b)(u)=bu$
for $b\in B,u\in \Sigma$, is an isomorphism, and $\mathbf{can}$ is a
surjective map;
\item $_B\Sigma$ is completely faithful, and
$\mathbf{can}$ is a bijective map.
\end{enumerate}
\end{theorem}

\begin{proof}
First we will prove the first statement. It is obvious that the
condition (1) implies the other conditions.

$(2)\Rightarrow(1)$ That $_A\coring{C}$ is flat follows from
$_B\Sigma$ is flat and $\mathbf{can}$ is an isomorphism of
$A$-bimodules. To prove that $(F,G)$ is a pair of equivalences, it
is enough to use the commutativity of the diagrams \eqref{counitmod}
and \eqref{unitmod}.

$(3)\Rightarrow(1)$ Follows from the Gabriel-Popescu Theorem
\cite[Theorem 3.7.9]{Popescu:1973}, and the commutativity of the
diagram \eqref{unitmod}.

$(4)\Rightarrow(1)$ From $\mathbf{can}$ is an isomorphism, it
follows that $-\cotensor{\coring{C}}(\Sigma^*\tensor{B}\Sigma)\simeq
-\cotensor{\coring{C}}\coring{C}$, and the following diagram is
commutative and each of its morphisms is an isomorphism:
$$
\xymatrix{N\tensor{B}(\Sigma\cotensor{\coring{C}}
\Sigma^*)\tensor{B}\Sigma
\ar[rr]^{\psi_N\tensor{B}\Sigma} \ar[d]^\simeq &&
\big((N\tensor{B}\Sigma)\cotensor{\coring{C}}\Sigma^*\big)
\tensor{B}\Sigma
\ar[d]^\simeq\\
N\tensor{B}\big(\Sigma\cotensor{\coring{C}}(\Sigma^*\tensor{B}
\Sigma)\big)\ar[rr]^\simeq
&&
(N\tensor{B}\Sigma)\cotensor{\coring{C}}(\Sigma^*\tensor{B}\Sigma).}
$$
Then $\psi_N$ is an isomorphism (since $-\tensor{B}\Sigma$ is
faithful). From the commutativity of the diagram
\eqref{unitcounitmod}, $\mathbf{can}$ is an isomorphism, and
$-\tensor{B}\Sigma$ is faithful, we have $\upsilon$ is an
isomorphism. Finally from the commutativity of the diagram
\eqref{unitmod}, $\upsilon_N$ is an isomorphism for every
$N\in\rmod{B}$.

Now we will prove the second statement. Obviously the condition (1)
implies the condition (4).

$(2)\Rightarrow(1)$ By \cite[Proposition 2.7]{Zarouali:2004},
$(\Sigma^*\tensor{B}-,\Sigma\cotensor{\coring{C}}-)$ is an adjoint
pair, and furthermore, $\Sigma_\coring{C}$ is generator if and only
if $\Sigma\cotensor{\coring{C}}-$ is faithful. From the
commutativity of the diagram \eqref{unitcounitmod}, and
$\Sigma\cotensor{\coring{C}}-$ is faithful, it follows that can is
an isomorphism. Hence (2) follows.

$(3)\Rightarrow(2)$ Since $_A\coring{C}$ is projective, and
$\mathbf{can}$ is surjective, $\mathbf{can}$ is a retraction in
$\lmod{A}$. Since $\coring{C}$ is coseparable, it follows that
$\mathbf{can}$ is a retraction in $\lcomod{\coring{C}}$. Then for
every $M\in \rcomod{\coring{C}}$,
$M\cotensor{\coring{C}}\mathbf{can}$ is a retraction in $\rmod{k}$,
and from the commutativity of the diagram \eqref{counitmod},
$\zeta_M$ is surjective for every $M\in \rcomod{\coring{C}}$. By
\cite[Theorem IV.3.1]{MacLane:1998}, $\Sigma_\coring{C}$ is a
generator.

$(4)\Rightarrow(3)$ We have that
$\Sigma\cotensor{\coring{C}}\mathbf{can}$ is bijective. From the
commutativity of the diagram \eqref{unitcounitmod},
$\upsilon\tensor{B}\Sigma$ is bijective. Since $_B\Sigma$ is
completely faithful, $\upsilon$ and $\phi$ are also bijective.
\end{proof}

\begin{remark}
In \cite[Proposition 5.6]{Caenepeel/DeGroot/Vercruysse:unp}, the
authors have a similar version of our second statement. They state
that if $\coring{C}_A$ is projective,
$B=\operatorname{End}_\coring{C}(\Sigma)$, and $\coring{C}$ is
coseparable, then $(F,G)$ defined as above is a pair of inverse
equivalences.
\end{remark}

In order to give a generalization of \cite[Theorem
3.5]{Takeuchi:1977} and \cite[Corollary 7.6]{AlTakhman:2002}, we
need the following result.

First we will recall from \cite{Zarouali:2004} the definition of a
coring having a duality. For details we refer to
\cite{Zarouali:2004}. If $\coring{C}_A$ is flat and $M \in
\rcomod{\coring{C}}$ is finitely presented as a right $A$-module,
then \cite[19.19]{Brzezinski/Wisbauer:2003} the dual left $A$-module
$\rdual{M} = \hom{A}{M}{A}$ has a left $\coring{C}$-comodule
structure
$$\rdual{M} \simeq
\hom{\coring{C}}{M}{\coring{C}} \subseteq \hom{A}{M}{\coring{C}}
\simeq \coring{C} \tensor{A} \rdual{M}.$$ ($f\mapsto
(f\tensor{A}\coring{C})\circ\rho_M$.) Now, if ${}_A\rdual{M}$ turns
out to be finitely presented and ${}_A\coring{C}$ is flat, then
$\ldual{(\rdual{M})} = \hom{A}{\rdual{M}}{A}$ is a right
$\coring{C}$-comodule and the canonical map $\sigma_M : M
\rightarrow \ldual{(\rdual{M})}$ is a homomorphism in
$\rcomod{\coring{C}}$. This construction leads to a duality
\[
\rDual: \rcomod{\coring{C}}_0 \leftrightarrows
\lcomod{\coring{C}}_0: \lDual
\]
between the full subcategories $\rcomod{\coring{C}}_0$ and
$\lcomod{\coring{C}}_0$ of $\rcomod{\coring{C}}$ and
$\lcomod{\coring{C}}$ whose objects are the comodules which are
finitely generated and projective over $A$ on the corresponding side
(this holds even without flatness assumptions of $\coring{C}$). Call
it the \emph{basic duality}.

\begin{proposition}\cite[Proposition 3.5]{Zarouali:2004}\label{lnoethrightadj}
Let $\coring{C}$ be an $A$-coring such that ${}_A\coring{C}$ and
$\coring{C}_A$ are flat. Assume that $\rcomod{\coring{C}}$ and
$\lcomod{\coring{C}}$ are locally noetherian categories. If
${}_A\rdual{M}$ and $\ldual{N}_A$ are finitely generated modules for
every $M \in \rcomod{\coring{C}}_f$ and $N \in
\lcomod{\coring{C}}_f$, then the basic duality extends to a right
adjoint pair $ \rDual: \rcomod{\coring{C}}_f \leftrightarrows
\lcomod{\coring{C}}_f: \lDual$. (see \cite{Colby/Fuller:1983} for
the definition of a right adjoint pair.)
\end{proposition}

Let $\coring{C}$ be a coring over $A$ satisfying the assumptions of
Proposition \ref{lnoethrightadj}. We will say that $\coring{C}$
\emph{has a duality} if the basic duality extends to a duality
$$\rDual: \rcomod{\coring{C}}_f \leftrightarrows
\lcomod{\coring{C}}_f: \lDual,$$ where for a Grothendieck category
$\cat{C}$, $\cat{C}_f$ stands for the full subcategory whose objects
are the finitely generated objects.

For instance, A cosemisimple coring, and a coring $\coring{C}$ over
a QF ring $A$ such that $_A\coring{C}$ and $\coring{C}_A$ are
projective, are corings having a duality.

\begin{lemma}\label{injectivecogenerator}
Let $N\in{}^\coring{C}\mathcal{M}^\coring{D}$ be a bicomodule.
Suppose that $A$ is a QF ring. \begin{enumerate}[(a)] \item If $N$
is an injector-cogenerator as a right $\coring{D}$-comodule, then
$N$ is an injective cogenerator in $\mathcal{M}^\coring{D}$.
\item If $\coring{D}$ has a duality, and $N$
is an injective cogenerator in $\mathcal{M}%
^\coring{D}$ such that $N_B$ is flat, then $N$ is an
injector-cogenerator as a right $\coring{D}$-comodule.
\end{enumerate}
\end{lemma}

\begin{proof}
(a) Since $A$ is a QF ring, then $A_A$ is an injective cogenerator.
Hence $N_\coring{D}\simeq (A\tensor{A}N)_\coring{D}$ is an injective
cogenerator.

(b) Let $X_A$ be an injective cogenerator module. Since $A$ is a QF
ring, $X_A$ is projective. We have then the natural isomorphism
\[
(X\otimes_AN)\square_{\coring{D}}-\simeq
X\otimes_A(N\square_{\coring{D}}-)
:{}^{\coring{D}}\mathcal{M}\rightarrow \mathcal{M}_{k}.
\]
By \cite[Proposition 3.8]{Zarouali:2004}, $N_{\coring{D}}$ and $X_A$
are faithfully coflat, and then $X\otimes_AN$ is faithfully coflat.
Once again by \cite[Proposition 3.8]{Zarouali:2004} ($X\otimes_AN$
is a flat right $B$-module), $X\otimes_AN$ is injective cogenerator
in $\mathcal{M}^\coring{D}$.
\end{proof}

\begin{theorem}\label{equivalence6}
Suppose that $_A\coring{C}$, $\coring{C}_A$, $_B\coring{D}$ and
$\coring{D}_B$ are flat, and let
$X\in{}^{\coring{C}}\mathcal{M}%
^{\coring{D}}$ and
$\Lambda\in{}^{\coring{D}}\mathcal{M}^\coring{C}$. If $\coring{C}$
and $\coring{D}$ have a duality, then the following statements are
equivalent:
\begin{enumerate}[(1)] \item $(-\square
_{\coring{C}}X,-\square_{\coring{D}}\Lambda)$ is a pair of inverse
equivalences with $X_B$ and $\Lambda_A$ are flat;
\item $(\Lambda
\cotensor{\coring{C}}-,X \cotensor{\coring{D}}-)$ is a pair of
inverse equivalences with $_AX$ and $_B\Lambda$ are flat;

If in particular $A$ and $B$ are QF rings, then (1) and (2) are
equivalent to \item $\Lambda$ is quasi-finite injective cogenerator
as a right
$\coring{C}%
$-comodule with $X_B$ and $\Lambda_A$ are flat,
$\e{\coring{C}}{\Lambda}\simeq \coring{D}$ as corings and $X\simeq
\cohom{\coring{C}}{\Lambda}{\coring{C}}$ in
$^{\coring{C}}\mathcal{M}^\coring{D}$;\item $X$ is quasi-finite
injective cogenerator as a right $\coring{D}%
$-comodule with $X_B$ and $\Lambda_A$ are flat,
$\e{\coring{D}}{X}\simeq \coring{C}$ as corings and $\Lambda\simeq
\cohom{\coring{D}}{X}{\coring{D}}$ in
$^{\coring{D}}\mathcal{M}^\coring{C}$;\item $\Lambda$ is
quasi-finite injective cogenerator on both sides,
$\e{\coring{C}}{\Lambda}\simeq \coring{D}$ as corings, and $X\simeq
\cohom{\coring{C}}{\Lambda}{\coring{C}}$ in
$^{\coring{C}}\mathcal{M}^\coring{D}$;\item $X$ is quasi-finite
injective cogenerator on both sides, $\e{\coring{D}}{X}\simeq
\coring{C}$ as corings, and $\Lambda\simeq
\cohom{\coring{D}}{X}{\coring{D}}$ in
$^{\coring{D}}\mathcal{M}^\coring{C}.$
\end{enumerate}
\end{theorem}

\begin{proof}
The equivalence between (1) and (2) follows from Theorem
\ref{equivalence4}, and the fact that if $(-\square
_{\coring{C}}X,-\square_{\coring{D}}\Lambda)$ is a Frobenius pair,
$X_{\coring{D}}$ and $\Lambda_{\coring{C}}$ are coflat if and only
if $X_B$ and $\Lambda_A$ are flat (see the proof of \cite[Theorem
3.11]{Zarouali:2004}).

$(1)\Leftrightarrow(3)$ Obvious from Proposition
\ref{equivalence3}$(I)$, and Lemma \ref{injectivecogenerator}.

$(1)\Leftrightarrow(4)$ The proof is analogous to that of
``$(1)\Leftrightarrow(3)$''.

$(5)\Rightarrow(3)$ Trivial.

 $(1)\Rightarrow(5)$ Obvious from the
equivalence of (1), (2) and (3).

$(1)\Leftrightarrow(6)$ Follows by symmetry.

\end{proof}

\begin{remark}
The second part of \cite[Corollary 7.6]{AlTakhman:2002} (and also
the second part of \cite[12.14]{Brzezinski/Wisbauer:2003}) is true
in a more general context. Let $_A\coring{C}$ be flat. Consider the
statements: \begin{enumerate}[(1)] \item
$\Lambda\in{}_B\mathcal{M}^\coring{C}$ be quasi-finite as a right
$\coring{C}$-comodule; \item $\hom{\coring{C}}{M_{0}}{\Lambda}$ is a
finitely generated left $B$-module, for every finitely generated
comodule $M_{0}\in\mathcal{M}^\coring{C}$; \item
$\hom{\coring{C}}{M_{0}}{\Lambda}$ is a finitely generated
projective left $B$-module, for every finitely generated comodule
$M_{0}\in\mathcal{M}^\coring{C}$. \end{enumerate} We have,
$(1)\Rightarrow(2)$ holds if $B$ is a QF ring, and the category
$\mathcal{M}^\coring{C}$ is locally finitely generated, and the
converse implication holds if in particular $B$ is a semisimple
ring, and the category $\mathcal{M}^\coring{C}$ is locally finitely
generated. $(1)\Leftrightarrow(3)$ holds if $\coring{C}$ is a
cosemisimple coring. Moreover, for the two cases, if (1) holds, then
for every comodule $M\in\mathcal{M}^\coring{C},$
$$\cohom{\coring{C}}{\Lambda}{M}\simeq
\underset{\underset{I}{\longrightarrow}}{\lim}\hom{\coring{C}}{M_i}{\Lambda}^*,$$
where $(M_i)_{i\in I}$ is the family of all finitely generated
subcomodules of $M$.

We will prove at the same time $(1)\Rightarrow(2)$ and
$(1)\Rightarrow(3)$. Let $M_{0}\in\mathcal{M}^\coring{C}$ be a
finitely generated comodule. We have,
$\cohom{\coring{C}}{\Lambda}{M_{0}}^*=\hom{B}{\cohom{\coring{C}}{\Lambda}{M_{0}}}{B}\simeq
\hom{\coring{C}}{M_{0}}{\Lambda}$. Since the functor
$-\tensor{B}{\Lambda}$ is exact and preserves coproducts, the cohom
functor $\cohom{\coring{C}}{\Lambda}{-}$ preserves finitely
generated (resp. finitely generated projective) objects. In
particular, $\cohom{\coring{C}}{\Lambda}{M_{0}}$ is a finitely
generated (resp. finitely generated projective) right $B$-module,
and $\hom{\coring{C}}{M_{0}}{\Lambda}$ so is. Therefore,
$\cohom{\coring{C}}{\Lambda}{M_{0}}\simeq\cohom{\coring{C}}
{\Lambda}{M_{0}}^{**}\simeq
\hom{\coring{C}}{M_{0}}{\Lambda}^*$. Hence,
$\cohom{\coring{C}}{\Lambda}{M}\simeq
\underset{\underset{I}{\longrightarrow}}{\lim}\cohom{\coring{C}}
{\Lambda}{M_i}\simeq
\underset{\underset{I}{\longrightarrow}}{\lim}\hom{\coring{C}}
{M_i}{\Lambda}^*$,
where $(M_i)_{i\in I}$ is the family of all finitely generated
subcomodules of $M$ (since the cohom functor preserves inductive
limits).

 The proof of the implication $(i)\Rightarrow(ii)$ of
\cite[Proposition 1.3]{Takeuchi:1977}, remains valid to prove
$(2)\Rightarrow(1)$ and $(3)\Rightarrow(1)$.
\end{remark}

\section{Applications to induction functors}\label{induction}

In this section we particularize our results in the previous section
to induction functors introduced in \cite{Gomez:2002}.

\medskip

A coring homomorphism \cite{Gomez:2002} from the coring
$\coring{C}$\ into the coring $\coring{D}$ is a pair
$(\varphi,\rho)$, where $\rho:A\rightarrow B$ is a homomorphism of
$k$-algebras and $\varphi:\coring{C}\rightarrow\coring{D}$ is a
homomorphism of $A$-bimodules such that
\[
\epsilon_{\coring{D}}\circ\varphi=\rho\circ
\epsilon_{\coring{C}}\qquad\textrm{and
}\qquad\Delta_{\coring{D}}\circ\varphi=
\omega_{\coring{D},\coring{D}}\circ(
\varphi\otimes_A\varphi)\circ\Delta_{\coring{C}},
\]
where
$\omega_{\coring{D},\coring{D}}:\coring{D}\otimes_A\coring{D}%
\rightarrow\coring{D}\otimes_B\coring{D}$ is the canonical map
induced by $\rho:A\rightarrow B.$

\medskip

Now we will characterize when the induction functor $ -\tensor{A} B
: \rcomod{\coring{C}} \rightarrow \rcomod{\coring{D}}$ defined in
\cite[Proposition 5.3]{Gomez:2002} is an equivalence of categories.
The right $\coring{D}$-comodule structure on $M \tensor{A} B$ is
defined by (using Sweedler's sigma notation)
\[
\rho_{M \tensor{A} B}(m \tensor{A}b) = \sum m_{(0)} \tensor{A} 1_B
\tensor{B} \varphi(m_{(1)})b,
\]
where $\rho_M(m) = \sum m_{(0)} \tensor{A} m_{(1)}$ is the coaction
of a right $\coring{C}$-comodule $M$.
\\ We also define the functor
$-\square_{\coring{D}}(B\tensor{A}\coring{C}):
\rcomod{\coring{D}}\rightarrow \rcomod{\coring{C}}$, where the left
comultiplicationon the left $B$-module $B\tensor{A}\coring{C}$ is
given by:
$$\lambda_{B\tensor{A}\coring{C}}:B\tensor{A}\coring{C}\rightarrow
\coring{D}\tensor{B}B\tensor{A}\coring{C}\simeq
\coring{D}\tensor{A}\coring{C},\quad b\tensor{A}c\mapsto \sum
b\varphi(c_{(1)})\tensor{A}c_{(2)},$$ where
$\Delta_{\coring{C}}(c)=\sum c_{(1)}\tensor{A}c_{(2)}$.
\\Moreover, if $_A\coring{C}$ is flat, then we have an
adjunction
$(-\otimes_AB,-\square_{\coring{D}}(B\tensor{A}\coring{C}))$ (see
\cite[24.11]{Brzezinski/Wisbauer:2003}).

\begin{theorem}\label{equivmor}
Let $(\varphi,\rho) :\coring{C}\rightarrow\coring{D}$ be a
homomorphism of corings such that $_A\coring{C}$ and $_B\coring{D} $
are flat. If $\coring{C}_A$ and $\coring{D}_B$ are flat (resp.
$\coring{C}$ and $\coring{D}$ are coseparable), then the following
statements are equivalent
\begin{enumerate}[(a)] \item
$(-\otimes_AB,-\square_{\coring{D}}(B\tensor{A}\coring{C}))$ is a
pair of inverse equivalences;\item the functor $-\otimes_AB$ is
exact and faithful (resp. faithful), and there exists an isomorphism
of $B$-corings $\coring{D}\simeq B\coring{C}B.$
\end{enumerate}
\end{theorem}

\begin{proof}
From the proof of \cite[Theorem 4.1]{Zarouali:2004},
$-\otimes_AB\simeq -\cotensor{\coring{C}}{(\coring{C}\otimes_AB)}$.
The use of Proposition \ref{equivalence3}, Proposition
\ref{comatrixproperty1}, and \cite[23.9]{Brzezinski/Wisbauer:2003}
achieves the proof.

\end{proof}

\begin{corollary}
Let $\coring{C}$ be an $A$-coring. Then the following statements are
equivalent
\begin{enumerate}[(a)]
\item The forgetful functor
$U_r:\mathcal{M}^{\coring{C}%
}\rightarrow\mathcal{M}_A$ is an equivalence of categories;
\item The forgetful functor
$U_l:{}^{\coring{C}}\mathcal{M}\rightarrow{}_A\mathcal{M}$ is an
equivalence of categories;\item there exists an isomorphism of
$A$-corings $A\simeq \coring{C}.$
\end{enumerate}
\end{corollary}

\begin{proof}
It is enough to apply the last theorem to the particular
homomorphism of corings
$(\epsilon_{\coring{C}},1_A):\coring{C}\rightarrow A$ which gives
the well-known adjunction $(U_r,-\tensor{A}{\coring{C}})$, and
observing that the map $A\coring{C}A\rightarrow \coring{C}$ defined
by $[a\otimes c\otimes a'\mapsto aca']$, is an isomorphism of
corings.
\end{proof}

Finally, given a homomorphism of corings, we give sufficient
conditions to have that the right induction functor is an
equivalence if and only if the left induction functor so is. Note
that for the case of coalgebras over fields (by (b)), or of rings
(well-known) (by (d)), we have the left right symmetry.

\begin{proposition}\label{equivmorsym}
Let $(\varphi,\rho):\coring{C}\rightarrow\coring{D}$ be a
homomorphism of corings such that $_A\coring{C}$, $_B\coring{D}$,
$\coring{C}_A$ and $\coring{D}_B$ are flat. Assume that at least one
of the following holds
\begin{enumerate}[(a)] \item $\coring{C}$ and
$\coring{D}$ have a duality, and $_AB$ and $B_A$ are flat;\item $A$
and $B$ are von Neumann regular rings;\item $B\otimes_A\coring{C}$
is coflat in
$^{\coring{D}%
}\mathcal{M}$ and $\coring{C}\otimes_AB$ is coflat in
$\mathcal{M}%
^{\coring{D}}$ and $_AB $ and $B_A$ are flat;\item $\coring{C}$ and
$\coring{D}$ are coseparable corings.
\end{enumerate} Then the following statements are
equivalent
\begin{enumerate}
\item $-\otimes_AB:\mathcal{M}^{\coring{C}%
}\rightarrow\mathcal{M}^\coring{D}$ is an equivalence of
categories;\item
$B\otimes_A-:{}^{\coring{C}}\mathcal{M}\rightarrow{}^{\coring{D}%
}\mathcal{M}$ is an equivalence of categories.
\end{enumerate}
\end{proposition}

\begin{proof}
Obvious from Theorem \ref{equivalence4}, Theorem \ref{equivalence5}
and Theorem \ref{equivalence6}.
\end{proof}

\section{Applications to entwined modules and graded
ring theory}\label{entwined}

In this section we particularize our results to corings associated
to entwined structures and in particular those associated to a
$G$-graded algebra and a right $G$-set, where $G$ is group. As in
\cite{Zarouali:2004}, we adopt the notations of \cite{DelRio:1992}
and \cite{Caenepeel/Militaru/Zhu:2002}.

\medskip

Let $(\alpha,\gamma):(A,C,\psi)\rightarrow(A',C',\psi')$ be a
morphism in $\mathbb{E}_\bullet^\bullet(k)$ (see
\cite{Caenepeel/Militaru/Zhu:2002}). We recall from
\cite{Zarouali:2004} that
$(\alpha\otimes\gamma,\alpha):A\otimes{C}\rightarrow{A'}\otimes{C'}$
is a morphism of corings and the functor $F$ defined in \cite[Lemma
8]{Caenepeel/Militaru/Zhu:2002} satisfies the commutativity of the
diagram
\[
\xymatrix{\mathcal{M}^{A\otimes
C}\ar[rr]^{-\otimes_AA'}\ar[d]^\simeq & &
\mathcal{M}^{A'\otimes C'}\ar[d]^\simeq \\
\mathcal{M}(\psi)_A^{C}\ar[rr]_{-\otimes_AA'} & &
\mathcal{M}(\psi')_{A'}^{C'},}
\]
where $-\otimes_AA':\mathcal{M}^{A\otimes C}\rightarrow
\mathcal{M}^{A'\otimes C'}$ is the induction functor defined in
\cite[Proposition 5.3]{Gomez:2002}.

We obtain the following result concerning the category of entwined
modules.

\begin{theorem}\label{entwined modules}
Let $(\alpha,\gamma):(A,C,\psi)\rightarrow(A',C',\psi')$ be a
morphism in $\mathbb{E}_\bullet^\bullet(k)$, such that $_{k}C$ and
$_{k}D$ are flat. If either $(A\otimes C)_A$ and $(A'\otimes
C')_{A'}$ are flat (e.g., if $\psi$ and $\psi'$ are isomorphisms),
(resp. $A\otimes C$ and $A'\otimes C'$ are coseparable (see
\cite[Theorem 38(1)]{Caenepeel/Militaru/Zhu:2002})), then the
following statements are equivalent
\begin{enumerate}[(a)] \item The functor
$-\otimes_AA':\mathcal{M}(\psi)_A^C\rightarrow
\mathcal{M}(\psi')_{A'}^{C'}$ defined in \cite[Lemma
8]{Caenepeel/Militaru/Zhu:2002} is an equivalence;
\item the functor
$-\otimes_AA'$ is exact and faithful (resp. faithful), and there
exists an isomorphism of $A'$-corings $A'\otimes C'\simeq
A'(A\otimes C)A'.$
\end{enumerate}
\end{theorem}

\begin{proof}
Follows immediately from Theorem \ref{equivmor}.
\end{proof}

\medskip

Now, let $G$ and $G'$ be two groups, $A$ be a $G$-graded
$k$-algebra, $A'$ be a $G'$-graded $k$-algebra, $X$ be a right
$G$-set, and $X'$ be a right $G'$-set. Let $\psi:kX\otimes
A\rightarrow A\otimes kX$ be the map defined by $[x\otimes
a_{g}\mapsto a_{g}\otimes xg]$. Analogously we define the map
$\psi':kX'\otimes A'\rightarrow A'\otimes kX'$. Let $kG$, $kG'$ be
the canonical Hopf algebras, and $kX$, $kX'$ be the canonical
grouplike coalgebras (see \cite{Caenepeel/Militaru/Zhu:2002}). From
\cite{Caenepeel/Militaru/Zhu:2002}, we have
$(kG,A,kX)\in\mathbb{DK}_\bullet^\bullet(k),$ $(A,kX,\psi)
\in\mathbb{E}_\bullet^\bullet(k) ,$ and
$\mathcal{M}(kG)_A^{kX}\simeq gr-( A,X,G).$  The comultiplication
and the counit maps of the coring $A\otimes kX$ are defined by:
$$\Delta_{A\otimes kX}(a\otimes x)=(a\otimes
x)\tensor{A}(1_A\otimes x), \quad \epsilon_{A\otimes kX}(a\otimes
x)=a \quad (a\in A,x\in X).$$

We recall that the corings $A\otimes kX$ and $A'\otimes kX'$ are
coseparable. A proof is obtained by using \cite[Corollary
3.6]{Brzezinski:2002} and \cite[Proposition
101]{Caenepeel/Militaru/Zhu:2002}. The proofs of the following
results are based on the results given in Section 6 of
\cite{Zarouali:2004}. We refer to \cite{DelRio:1992} for the
definitions of the terminology used in the following results.

Let $\widehat{A}=A\otimes kX$ be the $X\times X$-graded
$A-A$-bimodule associated to the $(A\otimes kX)-(A\otimes kX)$-
bicomodule $A\otimes kX$. It is clear that $\widehat{A}$ is
isomorphic as a bigraded bimodule to Del R\'{\i}o's
``$\widehat{A}$'' (see \cite{DelRio:1992}). The gradings are
$\widehat{A}_x=A\otimes kx$, and ${}_x\widehat{A}=\{\sum_ia_i\otimes
x_i\mid x_ig^{-1}=x, \forall i,\forall g\in G:(a_i)_g\neq 0\}$
($x\in X$).

\begin{corollary}\label{graded1}(\cite[Proposition
2.1]{DelRio:1992}) The following statements are equivalent
\begin{enumerate}[(1)] \item the categories
$gr-(A,X,G)$ and $gr-(A',X',G')$ are equivalent;
\item there are an $X\times X'$-graded $A-A'$-bimodule
$P$ and an $X'\times X$-graded $A'-A$-bimodule $Q$ such that
$$P\widehat{\otimes}_{A'}Q\simeq \widehat{A}\quad
\textrm{and}\quad
 Q\widehat{\otimes}_AP\simeq \widehat{A'}.$$
\end{enumerate}
Moreover, if $P$ and $Q$ satisfy the condition (2), then
$-\widehat{\otimes}_AP$ and $-\widehat{\otimes}_{A'}Q$ are inverse
equivalences.
\end{corollary}

\begin{proof}
Clear from Proposition \ref{equivalence1} (using the fact that the
corings $A\otimes kX$ and $A'\otimes kX'$ are coseparable) and
\cite[Corollary 6.5]{Zarouali:2004}.
\end{proof}

The following result is similar to \cite[Corollary
2.4]{DelRio:1992}.

\begin{corollary}
Let $A$ be a $G$-graded $k$-algebra, $X$ be a right $G$-set, and $B$
be an $k$-algebra. The following statements are equivalent
\begin{enumerate}[(1)] \item the category $gr-(A,X,G)$
is equivalent to $\rmod{B}$;
\item the category $(G,X,A)-gr$ is equivalent to
$\lmod{B}$;
\item there exists an $X_0\times X$-graded
$B-A$-bimodule $P$, such that $X_0$ is a singleton, $P$ is finitely
generated projective in $\rmod{A}$, and generator in $gr-(A,X,G)$,
and the morphism of $k$-algebras $\phi:B\to
\operatorname{End}_{gr-(A,X,G)}(P)$ defined by $\phi(b)(p)=bp$ for
$b\in B,p\in P$, is an isomorphism.
\end{enumerate}
\end{corollary}

\begin{proof}
Follows immediately from Theorem \ref{equivmodcomod}.
\end{proof}

In \cite[Theorem 2.3]{DelRio:1992}, A. Del R\'{\i}o gave a
characterization when $(-\widehat{\otimes}_AP,\h{P_{A'}}{-})$ is a
pair of inverse equivalences. We think that the following result
gives a simple characterization of it. Our result is also a
generalization of Morita's characterization of equivalence Remark
\ref{Morita}.

\begin{theorem}\label{graded2}
Let $P$ be an $X\times X'$-graded $A-A'$-bimodule. Then the
following are equivalent
\begin{enumerate}[(1)] \item
$(-\widehat{\otimes}_AP,\h{P_{A'}}{-})$ is a pair of inverse
equivalences;
\item \begin{enumerate}[(a)]\item $_xP$ is finitely
generated projective in $\rmod{A'}$ for every $x\in X$, and $P_{x'}$
is finitely generated projective in $\lmod{A}$ for every $x'\in X'$
(resp. and $P$ is a generator in both $gr-(A',X',G')$ and
$(G,X,A)-gr$),
\item the following bigraded bimodules maps:
 $\psi:\widehat{A}\rightarrow \h{P_{A'}}{P}$ defined
by $\psi(a\otimes x)(p)=a({}_xp)\in
 P$ $(a\in A,x\in X,p\in P)$, and
$\psi':\widehat{A'}\rightarrow \h{_AP}{P}$ defined by
 $\psi'(a'\otimes x')(p)=(p_{x'})a'\in
 P$ $(a'\in A',x'\in X',p\in P)$, are isomorphisms
(resp. are surjective maps).
\end{enumerate}

\item \begin{enumerate}[(a)]\item $_xP$ is finitely
generated projective in $\rmod{A'}$ for every $x\in X$, \item the
evaluation map
$$\varepsilon_{\widehat{A'}}:\h{P_{A'}}{\widehat{A'}}\widehat{\otimes}_AP\rightarrow
\widehat{A'},\; \varepsilon_{\widehat{A'}}(f\tensor{A}{P})=f(p)$$
($f\in \h{P_{A'}}{\widehat{A'}}_x, p\in {}_xP, x\in X$) is an
isomorphism,
\item the functor $-\widehat{\otimes}_AP$ is faithful.
\end{enumerate}
\end{enumerate}
\end{theorem}

\begin{proof}
From \cite[Proposition 1.2]{Menini:1993}, the unit and the counit of
the adjunction $(-\widehat{\otimes}_AP,\h{P_{A'}}{-})$ are given
respectively by $\eta_M:M\rightarrow
\h{P_{A'}}{M\widehat{\otimes}_AP},$ $\eta_M(m)(p)=\sum_{x\in
X}m_x\tensor{A}{_xp}$ ($m=\sum_{x\in X}m_x\in M,p=\sum_{x\in
X}{}_xp\in P$), and
$\varepsilon_N:\h{P_{A'}}{N}\widehat{\otimes}_AP\rightarrow N$,
$\varepsilon_N(f\tensor{A}{P})=f(p)$ ($f\in \h{P_{A'}}{N}_x,p\in
{}_xP, x\in X$). By \cite[Lemma 6.6(1)]{Zarouali:2004}, the functor
$\h{P_{A'}}{-})$ preserves inductive limits if and only if the
condition $(3)(a)$ holds. Then by \cite[Theorem 2.3]{Zarouali:2004}
(the coring $A'\otimes kX'$ is coseparable), there is a natural
isomorphism $\xymatrix{\delta:\h{P_{A'}}{-}\ar[r]^-\simeq &
-\widehat{\otimes}_{A'}\h{P_{A'}}{\widehat{A'}}}$. Set
$Q=\h{P_{A'}}{\widehat{A'}}$.

$(1)\Leftrightarrow(3)$ It follows from Proposition
\ref{equivalence3}(II) and Proposition \ref{comatrixproperty1}(2).

$(1)\Leftrightarrow(2)$ We can suppose that the condition $(2)(a)$
holds. We have $\eta_{\widehat{A}}:M\rightarrow
\h{P_{A'}}{\widehat{A}\widehat{\otimes}_AP},$
$\eta_{\widehat{A}}(a\otimes x)(p)=(a\otimes
x)\tensor{A}{}_xp=a(1_A\otimes
x)\tensor{A}{}_xp=a\lambda_P({}_xp)=\lambda_P(a({}_xp))=
\lambda_P(\psi(a\otimes
x)(p)),$ $(a\in A,x\in X,p\in P)$. Since $\delta$ is a natural
isomorphism, the following diagram is commutative:
$$
\xymatrix{\widehat{A}\ar[ddr]^{\psi}\ar[r]^-{\eta_{\widehat{A}}} &
\h{P_{A'}}{\widehat{A}\widehat{\otimes}_AP}
\ar[rr]^{\delta_{\widehat{A}\widehat{\otimes}_AP}}
& & \widehat{A}\widehat{\otimes}_AP\widehat{\otimes}_{A'}Q
\ar[r]^-\simeq & P\widehat{\otimes}_{A'}Q \\
& & & & \\
& \h{P_{A'}}{P}\ar[uu]^\simeq_{\h{P_{A'}}{\lambda_P}}
\ar[uurrr]^{\delta_P}& & & .}
$$
Therefore the unit of the adjunction
$(-\widehat{\otimes}_AP,-\widehat{\otimes}_{A'}Q)$ is
$$\xymatrix@1{1_{gr-(A,X,G)}\ar[r]^-\simeq &
-\widehat{\otimes}_A\widehat{A}\ar[rr]^-{-\widehat{\otimes}_A\psi_0}
& & -\widehat{\otimes}_AP\widehat{\otimes}_{A'}Q},$$ where $\psi_0$
is the map $\psi_0:\xymatrix{\widehat{A}\ar[r]^-{\psi} &
\h{P_{A'}}{P}\ar[r]^-{\delta_P} & P\widehat{\otimes}_{A'}Q}$. Hence
$-\widehat{\otimes}_AP$ is fully faithful if and only if $\psi$ is
an isomorphism. By \cite[Proposition 2.7]{Zarouali:2004},
$(Q\widehat{\otimes}_A-,P\widehat{\otimes}_{A'}-)$ is an adjoint
pair. Moreover $-\widehat{\otimes}_{A'}Q$ is fully faithful fully
faithfull" if and only if $P\widehat{\otimes}_{A'}-$ is fully
faithful, if and only if $\psi'$ is an isomorphism (see
\cite[Theorem IV.3.1]{MacLane:1998}).
\\ Finally, if the maps $\psi$ and $\psi'$ defined in
$(2)(b)$ are surjective, then, $P$ is a generator in $gr-(A',X',G')$
if and only if $-\widehat{\otimes}_{A'}Q$ is faithful, if and only
if $P\widehat{\otimes}_{A'}-$ is faithful, if and only if $\psi'$ is
an injective map (see \cite[Theorem IV.3.1]{MacLane:1998}). By
symmetry, $P$ is a generator in $(G,X,A)-gr$ if and only if $\psi$
is an injective map.
\end{proof}

Finally, let $f:G\rightarrow G'$ be a morphism of groups, $X$ be a
right $G$-set, $X'$ be a right $G'$-set, $\varphi:X\rightarrow X'$
be a map such that $\varphi(xg) =\varphi(x)f(g) $ for every $g\in
G,$ $x\in X.$ Let $A$ be a $G$-graded $k$-algebra, $A'$ be a
$G'$-graded $k$-algebra, and $\alpha:A\rightarrow A'$ be a morphism
of algebras such that $\alpha( A_{g}) \subset A_{f(g)}'$ for every
$g\in G.$

 We have, $\gamma:kX\rightarrow kX'$ such that
$\gamma(x) =\varphi(x)$ for each $x\in X$, is a morphism of
coalgebras, and $(\alpha ,\gamma) :(A,kX,\psi)
\rightarrow(A',kX',\psi') $ is a morphism in
$\mathbb{E}_\bullet^\bullet(k).$

 Let $T^*=-\otimes_AA': gr-( A,X,G)\rightarrow
gr-(A',X',G')$ be the functor defined in \cite[p.
531]{Menini/Nastasescu:1994}. $T^*$ makes commutative the following
diagram
\[%
\xymatrix{gr-(A,X,G)\ar[d]^\simeq \ar[rr]^{T^*} & &
gr-(A',X',G') \ar[d]^\simeq  \\
\mathcal{M}(kG)_A^{kX}
\ar[rr]^{-\otimes_AA'%
} & & \mathcal{M}(kG')_{A'%
}^{kX'}%
}
\]
(see Section 6 of \cite{Zarouali:2004}). Moreover, we have the
commutativity of the following diagram
\[%
\xymatrix{(G,X,A)-gr \ar[d]^\simeq \ar[rr]^{(T^*)'=A'\otimes_A-}
& & (G',X',A')-gr \ar[d]^\simeq\\
{}_A^{kX}\mathcal{M}(\psi^{-1})
\ar[d]^\simeq\ar[rr]^{A'%
\otimes_A-} & & {}_{A'}^{kX'}\mathcal{M}%
((\psi')^{-1}) \ar[d]^\simeq \\
{}^{kX\otimes A}\mathcal{M} \ar[d]^\simeq
\ar[rr]^{A'\otimes_A-}%
 & & {}^{kX'\otimes A'}\mathcal{M} \ar[d]^\simeq
\\
{}^{A\otimes kX}\mathcal{M} \ar[rr]^{A'\otimes_A-}%
 & & {}^{A'\otimes kX'}\mathcal{M}.}
\]

A. Del R\'{\i}o gave in \cite[Example 2.6]{DelRio:1992} an
interesting characterization when $T^*$ is an equivalence. Our
result gives an other characterization of it.

\begin{theorem}\label{graded3}
The following statements are equivalent\begin{enumerate} [(1)]
\item the functor $T^*:gr-(A,X,G)\rightarrow
gr-(A',X',G')$ is an equivalence; \item $T^*$ is faithful, and the
map
$$\omega:A'\tensor{A}\widehat{A}\tensor{A}A'\rightarrow
\widehat{A'}, \; a'\tensor{A}(a\otimes x)\tensor{A}a''\mapsto
a'\alpha(a)a''\otimes \varphi(x)g'$$ $(a\in A,a'\in A',x\in X, g'\in
G',a''\in A'_{g'})$, is bijective.
\end{enumerate}
\end{theorem}

\begin{proof}
Let $(\varphi,\rho) :\coring{C}\rightarrow\coring{D}$ be a
homomorphism of corings. From the proof of
\cite[24.11]{Brzezinski/Wisbauer:2003}, the counit of the adjunction
$(-\tensor{A}B,-\square_{\coring{D}}(B\tensor{A}\coring{C}))$ is
given by
$$\psi_N:(N\square_{\coring{D}}(B\tensor{A}\coring{C}))\tensor{A}B\rightarrow
N, \; \sum_in_i\tensor{A}c_i\tensor{A}b\mapsto
\sum_in_i\rho(\epsilon_{\coring{C}}(c_i))b$$ $(N\in
\rcomod{\coring{D}})$. This yields the map:
$$\omega:B\tensor{A}\coring{C}\tensor{A}B\rightarrow
\coring{D},\; b'\tensor{A}c\tensor{A}b\mapsto \sum
b'\varphi(c_{(1)})\rho(\epsilon_{\coring{C}}(c_{(2)}))b,\quad (c\in
\coring{C},b,b'\in B).$$ In our case, the last map is exactly that
mentioned in the condition (2). Finally, our result follows from
Theorem \ref {equivmodcomod}.
\end{proof}

Finally we give the following consequence of Proposition
\ref{equivmorsym}:

\begin{proposition}\label{graded4}
The following are equivalent\begin{enumerate} [(1)]\item the functor
$T^*:gr-(A,X,G) \rightarrow gr-( A',X',G')$ is an equivalence;
\item the functor $(T^*)':(G,X,A)-gr\rightarrow
(G',X',A')-gr$ is an equivalence.
\end{enumerate}
\end{proposition}

\section*{Acknowledgements}
I would like to thank my advisor, Professor Jos\'e
G\'omez-Torrecillas, for his time and advice.

\end{document}